\input amstex
\documentstyle{amsppt}
\magnification=\magstep1 \NoRunningHeads
\NoBlackBoxes

\input pictexwd.tex
\font\tiny=cmr5
\font\llarge=cmr42

\topmatter
\title On self-similarities of ergodic flows
\endtitle
\author
Alexandre I. Danilenko and Valery V. Ryzhikov
\endauthor

\thanks
The second named author was supported in part by the Program of Support of Leading Scientific Schools of RF (grant no.
8508.2010.1).
\endthanks

\abstract
Given an ergodic flow $T=(T_t)_{t\in\Bbb R}$, let $I(T)$ be the set of reals $s\ne 0$ for which the flows $(T_{st})_{t\in\Bbb R}$ and $T$ are isomorphic.
It is proved that  $I(T)$
is a Borel multiplicative subgroup of $\Bbb R^*$.
It carries a natural Polish group topology which is stronger than the topology induced from $\Bbb R$.
There exists a mixing flow $T$ such that $I(T)$ is an uncountable meager subset of $\Bbb R^*$.
 For a generic flow $T$, the transformations $T_{t_1}$ and $T_{t_2}$  are spectrally disjoint whenever $|t_1|\ne |t_2|$.
 A generic transformation (i) embeds  into a flow $T$ with $I(T)=\{1\}$ and (ii) does not embed  into a flow with $I(T)\ne \{1\}$.
For each countable multiplicative subgroup $S\subset\Bbb R^*$,
it is constructed a  Poisson  suspension flow $T$ with  simple spectrum such that $I(T)=S$.
If $S$ is without rational relations then there is
 a rank-one weakly mixing rigid flow $T$ with $I(T)=S$.
\endabstract

\address
 Institute for Low Temperature Physics
\& Engineering of National Academy of Sciences of Ukraine, 47 Lenin Ave.,
 Kharkov, 61164, UKRAINE
\endaddress
\email alexandre.danilenko\@gmail.com
\endemail

\address
 Department of Mechanics and Mathematics, Lomonosov
Moscow State University, GSP-1, Leninskie Gory, Moscow, 119991, Russian
Federation
\endaddress
\email vryzh\@mail.ru
\endemail


\endtopmatter
\document

\head 0. Introduction
\endhead

The isomorphism problem for measure preserving group actions is one of the central problems in ergodic theory.
Even within the framework of a single action this problem raises some interesting and difficult questions.
For instance, consider an action of $\Bbb R^n$, i.e. a multidimensional flow.
Then the automorphisms of $\Bbb R^n$  generate
(via linear time changes in the flow)
a continuum of new flows with possibly different classical invariants such as entropy and spectrum.
We also note that an $\Bbb R^n$-flow has a rich structure of subactions corresponding to lower dimensional subgroups and co-compact lattices.
A natural problem is to investigate (i) when these subactions are isomorphic, (ii) which invariants can distinguish non-isomorphic subactions, etc.
Our paper is devoted to the simplest particular case of this general problem.
We will only consider flows with one-dimensional time, i.e. $n=1$.
However even in this case there are a lot of open problems on inner symmetries and asymmetries of flows.

Let $T=(T_t)_{t\in\Bbb R}$ be an ergodic free measure preserving flow on a standard non-atomic probability space $(X,\goth B,\mu)$.
Given $s\in\Bbb R^*$, we denote by $T\circ s$ the  flow $(T_{st})_{t\in\Bbb R}$.
Let
$$
I(T):=\{s\in\Bbb R^*\mid T\circ s\text{ is isomorphic to }T\}.
$$
It is easy to see that  $I(T)$ is a multiplicative subgroup of $\Bbb R^*$.
If $I(T)\not\subset\{-1,1\}$ then $T$ is called {\it self-similar}.
We are also interested in a closely related invariant
$$
E(T):=\{t\in\Bbb R^*\mid T_t\text{ is isomorphic to }T_1\}.
$$
Of course, $I(T)\subset E(T)$.
In the present paper we investigate properties of the invariants $I(T)$ and $E(T)$.

The property of {\it total} self-similarity for flows, i.e. when $I(T)\supset\Bbb R_+^*$ is well known in ergodic theory.
If $T$ possesses this property  then  the maximal spectral type of $T$ is Lebesgue \cite{KaT, Proposition~1.23}.
In \cite{Ma}, it was shown that the total self-similarity implies mixing of all orders for the horocycle flow on any surface of constant negative curvature (in fact, for any flow acting by translations on a homogeneous space of a semisimple Lie group by a lattice).
A generalization of that result was obtained in \cite{Ry1}:
 given an arbitrary  ergodic flow $T$, if $E(T)$ has  positive Lebesgue measure then $T$ is mixing of all orders.
However it is unknown so far whether Leb$(E(T))>0$ implies that $I(T)\supset\Bbb R_+^*$.
Moreover, the following version of D.~Ornstein's question \cite{Th} is  open:
\roster
\item"---"
is there an  ergodic flow $T$ such that
 $E(T)=\Bbb R^*$ but $I(T)=\{1\}$?
\endroster
We note however that if $T$ acts by translations on a homogeneous space of a Lie group by a
lattice and $E(T)$ is uncountable then $I(T)\supset\Bbb R^*_+$ \cite{St}.

In \cite{dJR} it is proved that if $T$ is weakly mixing simple prime flow, then $I(T)=E(T)$ and $T\circ s\perp^{\text{F}} T$ whenever $s\not\in I(T)$.
The symbol $\perp^{\text F}$ denotes the disjointness in the sense of Furstenberg \cite{Fu}.
It was stated in \cite{Ry4} that for a rank-one mixing flow $T$ we have $T\circ s\perp^{\text{F}} T$ for all $s\in\Bbb R^*\setminus\{1\}$.
It is interesting to note that if $T$, in addition,  has a simple Lebesgue spectrum (see \cite{Pr} in this connection)
then  the  Koopman representations
of $T\circ s$ and $T$
are unitarily equivalent for all $s\in\Bbb R^*$.
This illustrates a drastic difference between  the disjointness in the sense of Furstenberg and the spectral disjointness.

In the case when the maximal spectral type of a flow $T$ is singular, the transformation $T_t$ is spectrally disjoint with $T_1$ for Leb.-a.a. $t\in \Bbb R$ \cite{Ry3}.
An example of a non-mixing flow  $T$
with minimal self-joinings and $I(T)=\{1\}$ was constructed in \cite{dJP} (see also \cite{dJR, Proposition~6.8}).

It is easy to construct a non-mixing  ergodic  flow $T$ with  $I(T)$ infinite.
 Consider, for instance, an infinite Cartesian product
$$
\Bbb R\ni t\mapsto   T_t:=\cdots\times S_{\alpha^{-1}t}\times S_{t}\times S_{\alpha t}\times S_{\alpha^2t}\times\cdots,
$$
where $S$ is a non-mixing weakly mixing flow.
Then $I(T)\supset\{\alpha^n\mid n\in\Bbb Z\}$.
In this connection a natural question was posed in \cite{Ry1}:
\roster
\item"---"
is there a non-mixing weakly mixing flow $T$ such that $E(T)$ is uncountable.
\endroster
It remains open (cf. Theorem~2.1 below).

An extensive study of various  self-similarity problems was undertaken in a recent paper \cite{FrL} by K.~Fr\c{a}czek and M.~Lema\'nczyk.
In particular, the following was done.
\roster
\item"(i)" Examples of non-self-similar ergodic special flows built over certain interval exchange transformations
are given.
They include some non-mixing smooth flows on translation surfaces (see also \cite{Ku} for constructions of smooth non-self-similar flows on each surface of genus $\ge 2$).
\item"(ii)"
For each countable subgroup $G\subset\Bbb R^*$, a weakly mixing flow $T$ is constructed with $I(T)=G$.
\item"(iii)" If $G$ is independent as a subset of the $\Bbb Q$-linear space $\Bbb R$ then  a weakly mixing Gaussian flow $T$ with simple spectrum
is constructed such that  $I(T)=G\sqcup(-G)$ and $T\circ t$ is spectrally disjoint with $T$
whenever $|t|\not\in G$.
\endroster
They also raised several questions as a certain  program for further investigation of self-similarity problems.
These questions together with  additional ones kindly sent to us
 by M.~Lema\'nczyk stimulated our present work.
We give a complete answer to the following.
\roster
\item"---"{(Q1)} Are the sets $I(T)$ and $E(T)$  Borel subsets of $\Bbb R^*$?
\item"---"{(Q2)} Is there a natural Polish topology on $I(T)$?
\item"---"{(Q3)} Is there a flow $T$ for which the group $I(T)$ is  uncountable and has zero Lebesgue measure?
\item"---"{(Q4)} Can we embed a typical  transformation into a self-similar flow?
\item"---"{(Q5)} Is the absence of self-similarity generic in the set of all measure preserving flows on $(X,\goth B,\mu)$?
\item"---"{(Q6)} Find  weakly mixing rank-one self-similar flows.
\endroster
 Thus we solved Problems 1, 2, 5, 6 from the list in \cite{FrL, Section~10} plus two  additional ones (Q2) and (Q6) by M.~Lema\'nczyk.
We also remove a redundant  ``independence'' condition on $G$ from (iii).
Moreover,  (iii) is proved in the full generality for both Poisson and Gaussian flows.
We now state precisely these and other main results of our work.

\proclaim{Proposition 1.3}
\roster
\item"(i)"
$I(T)$ and $E(T)$ are  Borel subsets of $\Bbb R^*$.
\item"(ii)"
There is a topology $\tau$ on $I(T)$ which is stronger than the topology induced from $\Bbb R^*$ and such that $(I(T),\tau)$ is a Polish topological group.
\endroster
\endproclaim

It is possible to have $I(T)\ne E(T)$.
 The set $E(T)$ does not need to be a subgroup of $\Bbb R^*$ for an arbitrary  $T$ (see Example~1.1 and  S.~Tikhonov's Example~1.2).

\proclaim{Theorem 2.1}
There is a mixing (of all orders) flow $T$ such that  $I(T)$
is uncountable but $I(T)\not\supset\Bbb R_+^*$.
\endproclaim

By Aut$(X,\mu)$ and Flow$(X,\mu)$
we denote the group of $\mu$-preserving transformations of $(X,\mu)$
and the set of $\mu$-preserving flows on $(X,\mu)$ respectively.
We endow these sets with the natural Polish topologies (see Section~1 below).
As usual, we say that a property is {\it generic} in a Polish space $P$ if the subset of elements satisfying this property is residual  in $P$.

\proclaim{Corollary 3.3}   For a generic flow $T$, the group $I(T)$ is trivial.
Moreover, $\sigma_T\perp\sigma_{T\circ t}$ if  $|t|\ne 1$ and
$T\perp^{\text F} T\circ t$ if $t=-1$.
\endproclaim

Here and below $\sigma_T$ denotes a measure of the maximal spectral type of $T$.

\proclaim{Theorem 3.6}
\roster
\item"(i)"
 A generic transformation from $\text{\rom{Aut}}(X,\mu)$ embeds into a flow $T$ such that  $I(T)=\{1\}$.
Moreover, it embeds into a flow possessing all the properties listed in Corollary~3.3.
\item"(ii)"
 A generic transformation from $\text{\rom{Aut}}(X,\mu)$ does not embed into a flow $T$ with $I(T)\not=\{1\}$.
\endroster
\endproclaim

We note that (ii) does not follow directly from (i)  because a  generic transformation from $\text{\rom{Aut}}(X,\mu)$ embeds into continuum of pairwise non-isomorphic flows~\cite{SE}.

\proclaim{Theorem 4.1}
Let $S$ be a countable subgroup of $\Bbb R_+^*$ such that $S$ considered as a subset of the $\Bbb Q$-linear space $\Bbb R$ is independent.
Denote by $T^S$ the Cartesian product flow $\bigotimes_{s\in S}T\circ s$ acting on the space $(X,\mu)^S$.
For a generic flow $T\in\text{\rom{Flow}}(X,\mu)$,
\roster
\item"\rom(i)" the flow
$T^S$ is rank one rigid and weakly mixing,
\item"\rom(ii)"
 $I(T^S)=S$ and, moreover,
\item"\rom(iii)"
$T^S\perp^{\text{\rom F}} (T^S)\circ t$ for each real $t\not\in S\cup\{0\}$.
\endroster
\endproclaim

For arbitrary countable subgroups of $\Bbb R^*$ we prove the following.

\proclaim{Theorem 4.3}
Let $S$ be a countable subgroup of $\Bbb R^*$.
There is a weakly mixing rank-one rigid flow $T$ such that $I(T)\supset S$.
\endproclaim

\proclaim{Theorems 4.4, 4.10}
Let $S$ be a countable subgroup of $\Bbb R^*$.
There is a weakly mixing Poisson suspension  flow $\widetilde W$ with a simple spectrum such that $I(\widetilde W)=S$ and $\sigma_{\widetilde W}\perp \sigma_{\widetilde W\circ t}$ for each positive $t\not\in S$.
Hence  there is also a weakly  mixing Gaussian flow $F$
with a simple spectrum such that
$I(F)=S\cup(-S)$ and $\sigma_F\perp\sigma_{F\circ t}$ for each  positive $t\not\in S$.
\endproclaim

It follows, in particular, that if $S$ is contained in $\Bbb R^*_+$
 then the corresponding Poisson suspension flow  is not   isomorphic to
its inverse.
This is in a strong contrast with the Gaussian case: $I(F)\ni -1$ for each Gaussian flow $F$.

\remark{Remark\footnote{Suggested by the referee.}}
Theorem~4.10 brings the answer to a long standing open question from harmonic analysis\footnote{See  \cite{Le}, Section ``Future directions".}: {\it there is a Gaussian system with a simple spectrum
such that the spectral measure $\sigma$ of the Gaussian process (which determines the Gaussian system) can not be concentrated on any subset
without rational relations}.
Indeed, if $S=\{2^n\mid n\in\Bbb Z\}$
and $\sigma(A)=1$ for a subset $A\subset\Bbb T$ then $\sigma(A\cap 2A)=1$.
\endremark

The authors thank M.~Lema\'nczyk for his questions and useful discussions and S.~Tikhonov who showed us  Example~1.2.
We are also grateful to J.-P.~Thouvenot and E.~Roy for their helpful comments.
We thank the anonymous referee for the valuable remarks and suggestions which improved the paper.

\head  1.  Topological and algebraic properties of $I(T)$ and $E(T)$
\endhead

Let $T=(T_t)_{t\in\Bbb R}$ be an ergodic free measure preserving flow on a standard probability space $(X,\goth B,\mu)$.
In this section we study  algebraic and topological properties of $I(T)$ and $E(T)$ and answer questions (Q1) and (Q2).

Denote by
 $\Lambda(T)\subset\Bbb R$  the discrete spectrum of $T$.
Then $\Lambda(T)$ is a countable subgroup of $\Bbb R$.
Denote by $\goth F$ the Kronecker factor of $T$,
i.e.  $\goth F\subset\goth B$ is the sub-$\sigma$-algebra  generated by all proper functions of $T$.
It is easy to verify that
$$
s\Lambda(T)=\Lambda(T)\qquad\text{for each}\quad s\in I(T).\tag1-1
$$

We first give a simple example of a free  ergodic  flow  $T$ such that $I(T)\ne E(T)$.
This flow has a non-trivial discrete spectrum.
Weakly mixing flows with this property also exist but they are more involved (see Example~1.2 below).

\example{Example 1.1}
Let $B=(B_t)_{t\in\Bbb R}$ be a Bernoulli flow with infinite entropy and let $P=(P_t)_{t\in\Bbb R}$ be an ergodic flow with pure point spectrum $\Bbb Z$.
Then the product flow $B\times P$ is free and ergodic.
Since $\Lambda (B\times P)=\Bbb Z$, it follows from \thetag{1-1} that
$I(B\times P)\subset\{-1,1\}$.
The converse inclusion is obvious.
Hence $I(B\times P)=\{-1,1\}$.
If  $0\ne t\in\Bbb Z$ then  $B_t\times P_t=B_t\times\text{Id}$.
The ergodic components of this transformation are all isomorphic to $B_t$. It follows  from the Ornstein's isomorphism theory for Bernoulli transformations that $E(B\times P)\supset\Bbb Z\setminus\{0\}$.
If $t\not\in\Bbb Z$ then the ergodic components of the transformation $B_t\times P_t$
have  non-trivial point spectrum.
Hence $t\not\in E(B\times P)$.
Thus $E(B\times P)=\Bbb Z\setminus\{0\}\ne I(B\times P)$.
\endexample

Denote by $\text{Aut}(X,\mu)$ the group of all $\mu$-preserving invertible transformations of $(X,\mu)$.
Endow it with the weak topology in which $R_n\to R$ if
$\mu(R_nA\triangle RA)\to 0$ as $n\to\infty$.
Then $\text{Aut}(X,\mu)$ is a Polish group \cite{Ha}.

In the following example by S.~Tikhonov,  a weakly mixing flow $V$ is constructed such that $I(V)\ne E(V)$ and $E(V)$ is not a subgroup of $\Bbb R^*$.

\example{Example 1.2}
It can be deduced easily from \cite{dRdS} and \cite{dJL} that there is a
residual subset $\Cal F\in\text{Aut}(X,\mu)$ such that for each transformation $S\in\Cal F$ the following holds.
\roster
\item"(i)" $S$ is weakly mixing.
\item"(ii)" There exists a flow $T\in\text{Flow}(X,\mu)$ such that
$T_1=S$,
\item"(iii)"
$C(S)=C(T_q)$  for each $q\in\Bbb Q$ and
\item"(iv)"
the centralizer of the infinite product transformation
$$
\cdots\times T_{\frac 14}\times T_{\frac 12}\times T_1\times T_{2}\times T_{4}\times\cdots
$$
of the product space $(X,\mu)^\Bbb Z$
is the infinite product $ C(S)^\Bbb Z$, i.e. this centralizer is as ``small'' as possible.
\endroster
According to \cite{Ti}, a generic transformation in Aut$(X,\mu)$ has a continuum of  roots in each residual subset of
Aut$(X,\mu)$.
Therefore there are  transformations $S\ne\widehat S$ in $\Cal F$ such that $S^2=\widehat S^2$.
Let $T$ be a flow satisfying (ii)--(iv)
and let $\widehat T$ be a flow satisfying (ii)--(iv) with $\widehat S$
instead of $S$.
We now define a flow $V$ on the space $(X,\mu)^\Bbb Z$ by setting
$$
V=
\cdots\times T\circ 2^{-2}\times T\circ 2^{-1}\times T\times  \widehat T\circ 2\times T\circ 2^2\times\cdots
$$
It follows from (i) that  $V$ is weakly mixing.
It is straightforward that $RV_1R^{-1}=V_2$, where $R:X^\Bbb Z\to X^\Bbb Z$ denotes the shift.
Therefore $2\in E(V)$.
We now show that $2^{-1}\not\in E(V)$.
Indeed, if $QV_1Q^{-1}=V_{2^{-1}}$ for some transformation $Q$ of $(X,\mu)^\Bbb Z$ then
$QV_2Q^{-1}=V_{1}$.
Hence the transformation $QR$ commutes with $V_1$ and
$$
R^{-1}V_1R=(QR)^{-1}V_{2^{-1}}QR.
$$
However then it follows from (iv) and (iii) that the transformations
$T_1$ and $\widehat T_1$, i.e. $S$ and $\widehat S$ in view of (ii), are conjugate by an element of the group $C(S)$.
Hence $S=\widehat S$, a contradiction.
Thus, $2^{-1}\not\in E(V)$ and therefore $2\notin I(V)$.
\endexample

Let $d_w$ be a complete metric on $\text{Aut}(X,\mu)$ compatible with the weak topology.
Denote by  $\text{Flow}(X,\mu)$ the set of all $\mu$-preserving flows on $(X,\mu)$.
Endow it with the topology of uniform weak convergence on the compact subsets in $\Bbb R$.
This topology is compatible with the following metric $d$:
$$
d(T,S)=\sup_{0\le t\le 1}d_w(T_t,S_t).
$$
Then $(\text{Flow}(X,\mu),d)$  is a Polish space (see \cite{dRdS}).

\proclaim{Proposition 1.3}
\roster
\item"(i)"
$I(T)$ and $E(T)$ are  Borel subsets of $\Bbb R$.
\item"(ii)"
There is a topology $\tau$ on $I(T)$ which is stronger than the topology induced from $\Bbb R^*$ and such that $(I(T),\tau)$ is a Polish topological group.
\endroster
\endproclaim
\demo{Proof}
There are two  commuting  actions of $\text{Aut}(X,\mu)$ and $\Bbb R^*$ on $\text{Flow}(X,\mu)$:
$$
\align
\text{Aut}(X,\mu)\times\text{Flow}(X,\mu)&\ni(R,T)\mapsto R\bullet T\in\text{Flow}(X,\mu) \text{ and} \tag1-2\\
\Bbb R^*\times\text{Flow}(X,\mu)&\ni(s,T)\mapsto T\circ s\in\text{Flow}(X,\mu),
\endalign
$$
where the flow $R\bullet T$ is given by $(R\bullet T)_t:=RT_tR^{-1}$ for all $t\in\Bbb R$.
The two actions are continuous.
We verify the continuity of the first one\footnote{By an advice of the referee.}.
Fix $R\in \text{Aut}(X,\mu)$, $T\in\text{Flow}(X,\mu)$ and $\epsilon>0$.
Since $(\text{Aut}(X,\mu), d_w)$ is a topological group and the mapping $[0,1]\ni t\mapsto T_t\in\text{Aut}(X,\mu)$ is uniformly continuous,
there exist reals $a_1,\dots,a_n\in[0,1]$ and $\delta>0$ such that
\roster
\item"$(\diamond)$"
$\sup_{0\le t\le 1}\min_{1\le i\le 1}d(T_t,T_{a_i})<\delta/2$ and
\item"$(\diamond)$"
if $d_w(R,\widetilde R)<\delta$ and $d_w(Q,T_{a_i})<\delta$ for some transformations $\widetilde R,Q\in\text{Aut}(X,\mu)$ and some $1\le i\le n$ then
$
d_w(\widetilde RQ\widetilde R^{-1}, RT_{a_i}R^{-1})<\epsilon
$
\endroster
Now take a flow $\widetilde T$ with $d(T,\widetilde T)<\delta/2$.
For each $t\in[0,1]$, we find $i$ such that $d_w(T_t,T_{a_i})<\delta/2$.
Then $d_w(\widetilde R\widetilde T_t\widetilde R^{-1},RT_{a_i}R^{-1})<\epsilon$ and hence
$
d(\widetilde R\bullet\widetilde T,R\bullet T)\le 2\epsilon.
$
Thus the action \thetag{1-2} is continuous.

It follows that the set
$$
A_T:=\{(R,s)\in\text{Aut}(X,\mu)\times\Bbb R^*\mid R\bullet T=T\circ s\}
$$
is a closed subgroup of $\text{Aut}(X,\mu)\times\Bbb R^*$.
The group $I(T)$ is the image of $A_T$ under a continuous homomorphism
$\pi:A_T\ni(R,s)\mapsto s\in\Bbb R^*$.
We may view $\pi$ as a one-to-one continuous homomorphism from the Polish group
$:A_T/\text{Ker\,}\pi$ onto $I(T)$.
  Thus $I(T)$ is a Borel subset of $\Bbb R$ since it is a one-to-one continuous image of the Polish space.
Let $\tau$ be the topology on $I(T)$ in which the map $:A_T/\text{Ker\,}\pi\to I(T)$
is a homeomorphism.
Then (ii) holds.

Let $M_T:=\{(R,s)\in \text{Aut}(X,\mu)\times\Bbb R^*\mid T_s=RT_1R^{-1}\}$.
Then $M_T$ is a closed subset of $\text{Aut}(X,\mu)\times\Bbb R^*$.
Since the centralizer $C(T_1)$ of  $T_1$, i.e. the  group of all
transformations commuting with $T_1$, is closed in $\text{Aut}(X,\mu)$,
there is a Borel subset $B\subset\text{Aut}(X,\mu)$ such that every
transformation $S\in\text{Aut}(X,\mu)$ can be written uniquely as a product $S=WR$ with $R\in B$ and $W\in C(T_1)$ \cite{Ke, 12.17}.
Now it is easy to verify that $E(T)$ is a one-to-one image of the Borel subset $M_T\cap(B\times \Bbb R^*)$ under a continuous map
$\tau:M_T\ni(R,s)\mapsto s\in\Bbb R^*$.
Hence $E(T)$ is Borel.
\qed
\enddemo

It follows from the Banach-Kuratowski-Pettis theorem \cite{Kel} that every Borel subgroup of $\Bbb R^*$ is either non-empty open or meager and of zero Lebesgue measure, we deduce from Proposition~1.3(ii) that

\proclaim{Corollary 1.4}
 Either $I(T)$ contains $\Bbb R^*_+$ and then the maximal spectral type of $T$ is Lebesgue (see \cite{KaT, Proposition~1.23}) or
 $I(T)$ is meager and  Leb$(I(T))=0$.
\endproclaim

\remark{Remark \rom{1.5}}
It follows from the proof of Proposition~1.3 that the following sequence of Polish groups
$$
1\to\text{Ker}\,\pi \overset\text{id}\to\longrightarrow M_T \overset\pi\to\longrightarrow I(T)\to 1\tag1-3
$$
is exact.
An interesting question is when it splits, i.e. there is a continuous
homomorphism $R:I(T)\ni s\mapsto R_s\in\text{Aut}(X,\mu)$ such that
$T_{st}=R_sT_tR_s^{-1}$ for all $t\in\Bbb R$ and $s\in I(T)$.
It splits in the  case when $T$ is a horocycle flow or when $T$ is a Bernoulli flow with infinite entropy.
Also, if $I(T)$ is isomorphic to $\Bbb Z^p$ then~\thetag{1-3} splits.
We do not know the answer in the general case.
\endremark

\head 2. Mixing flow $T$ with $I(T)$ meager and uncountable
\endhead

Our main purpose in this section is to  construct an ergodic flow $T$ such that $I(T)$ is uncountable and meager.
This answers  (Q3).
We construct such a flow as a 2-point extension of a horocycle  flow.
The extension is chosen in such a way to  partially ``destroy''  self-similarities of the base flow.
This means that uncountably many of elements of the corresponding geodesic flow in the base lift to the extension and some elements do not lift.
Measurable orbit theory plays a key role in choosing such an extension.

 We  first observe that if $T$ is ergodic and $I(T)$ is uncountable then $T$ is weakly mixing.
Indeed, this follows from \thetag{1-1}.

\proclaim{Theorem 2.1}
There is a mixing (of all orders) flow $T$ such that  $I(T)$
is uncountable but $I(T)\not\supset\Bbb R_+^*$.
\endproclaim

The proof of this theorem is based heavily on the orbit theory of amenable dynamical systems.
Therefore we begin this section with a preliminary material on the orbit theory.

If an equivalence relation $\Cal R$ on $(X,\goth B,\mu)$ is the orbit equivalence relation of a $\mu$-preserving action $T$ of a locally compact second countable group $G$ then $\Cal R$ is called {\it measure preserving}.
If every $\Cal R$-saturated measurable subset of $X$ is either $\mu$-null or $\mu$-conull then $\Cal R$ is called {\it ergodic}.
We note that $\Cal R$ is ergodic if and only if $T$ is ergodic.

If the $\Cal R$-class of a.e. point $x\in X$ is countable then $\Cal R$ is called {\it discrete}.
If the $\Cal R$-class of a.e. point $x\in X$ is uncountable then $\Cal R$ is called {\it continuous}.
If $\Cal R$ is ergodic then it is either discrete or continuous.

We do not give here the general definition of amenability  for equivalence relations (see \cite{Zi}) but just note that if $G$ is amenable then $\Cal R$ is {\it amenable}.
Given a  compact second countable group $K$, a Borel map $\alpha:\Cal R\to K$ is called a {\it cocycle} of $\Cal R$ if there is a $\mu$-conull subset $Y\subset X$ such that
$$
\alpha(x,y)\alpha(y,z)=\alpha(x,z)\quad\text{for all }(x,y),(y,z)\in\Cal R_Y:=\Cal R\cap (Y\times Y).
$$
We do not distinguish between cocycles which agree a.e.
Recall that two cocycles $\alpha,\beta:\Cal R\to K$ {\it agree a.e.}
if there is a $\mu$-conull  subset $Z\subset X$ such that
$\alpha(x,y)=\beta(x,y)$ for all $(x,y)\in\Cal R_Z$.
Two cocycles $\alpha,\beta:\Cal R\to K$ are {\it cohomologous} (we will denote $\alpha\approx\beta$) if there is a Borel map $\phi:X\to K$ and a $\mu$-conull subset $Y$ such that
$$
\alpha(x,x')=\phi(x)\beta(x,x')\phi(x')^{-1}\quad\text{for all }(x,x')\in\Cal R_Y.
$$
An invertible  $\mu$-preserving transformation $S$ of $X$ is called an {\it automorphism} of $\Cal R$ if there are $\mu$-conull subsets $X_1$ and $X_2$ such that $(S\times S)\Cal R_{X_1}=\Cal R_{X_2}$.
Then we can define a cocycle $\alpha\circ S:\Cal R\to K$ by setting
$\alpha\circ S(x,x'):=\alpha(Sx,Sx')$.
It is easy to verify that  $\alpha\approx\beta$ if and only if  $\alpha\circ S\approx\beta\circ S$.

A Borel measure preserving action $D=(D_h)_{h\in H}$ of a locally compact second countable group $H$ on $(X,\mu)$ is called {\it strictly $\Cal R$-outer} if there is a conull subset $X'\subset X$ such that
\roster
\item"(i)"
$D_h$ is an automorphism of $\Cal R$ for each $h\in H$ and
\item"(ii)"
if $(D_hx,x)\in\Cal R$ for some $x\in X'$ and $h\in H$ then $h=1_H$.
\endroster

A $\mu$-preserving invertible transformation $S$ of $X$ is called $\Cal R$-{\it inner}
if $(x,Sx)$ belongs to $\Cal R$ for a.a. $x$.
Of course, $S$ is an automorphism of $\Cal R$.
It is straightforward  that $\alpha\circ S\approx\alpha$ for each cocycle $\alpha$ of $\Cal R$ and  each $\Cal R$-inner automorphism $S$.
 Consider a measure preserving transformation $S^\alpha$ of the product space $(X\times K,\mu\times\lambda_K)$ by setting
$$
S^\alpha(x,k)=(Sx,\alpha(Sx,x)k),
$$
where $\lambda_K$ is the Haar measure on $K$.
Then $S^\alpha$ is called the $\alpha$-{\it skew product extension} of $S$.
If $\Cal R$ is generated by a  $G$-action $T=(T_g)_{g\in G}$  then $((T_g)^\alpha)_{g\in G}$ is a measure preserving $G$-action on   $(X\times K,\mu\times\lambda_K)$.
Denote by $\Cal R(\alpha)$ the orbit equivalence relation of this action.
It does not depend on a particular choice of $T$ generating $\Cal R$.
We note that $(x,k)\sim_{\Cal R(\alpha)}(x',k')$ if and only if $x\sim_\Cal R x'$ and $k'=\alpha(x',x)k$.
If $\Cal R(\alpha)$ is ergodic then $\alpha$ is called {\it ergodic}.

\proclaim{Proposition 2.2}
There are an amenable ergodic measure preserving continuous equivalence relation $\Cal T$ on a standard probability space $(\widetilde Y,\widetilde{\goth B},\widetilde\nu)$,
a strictly $\Cal T$-outer flow $\widetilde F=(\widetilde F_t)_{t\in\Bbb R}$ on $\widetilde Y$
and an ergodic cocycle $\widetilde\beta:\Cal T\to\Bbb Z/2\Bbb Z$ such that the set $L:=\{t\in\Bbb R\mid\widetilde\beta\circ\widetilde F_t\approx\widetilde\beta\}$ is a proper uncountable subgroup of $\Bbb R$.
\endproclaim
\demo{Proof}
We use a 3-step construction.

{\bf (A)}
Let $F'=(F'_t)_{t\in\Bbb R}$ denote the following flow on $(\Bbb T,\lambda_\Bbb T)$:
$$
F'_tz=z+t.
$$
We view $\Bbb T$ as the interval $[0,1)$.
The addition is considered  mod 1.
We note that
$F'$ is transitive and periodic, $F_{t+1}'=F_t'$ for each $t$.
Fix an irrational number $\theta_1$.
Denote by $\Cal R_{\theta_1}$ the orbit equivalence relation of the transformation $F_{\theta_1}'$ on $(\Bbb T,\lambda_\Bbb T)$.
Then $\Cal R_{\theta_1}$ is  discrete and ergodic.
There is a bijection between the  cocycles of $\Cal R_{\theta_1}$ with values in $\Bbb Z/2\Bbb Z$ and the set $\Cal M(\Bbb T,\Bbb Z/2\Bbb Z)$ of measurable functions from $\Bbb T$ to $\Bbb Z/2\Bbb Z$.
Such a bijection is defined in a highly non-unique way.
For instance, it is established by the map
$\beta\mapsto a_\beta$, where $a_\beta(z):=\beta(z,z+\theta_1)$.
The set $\Cal M(\Bbb T,\Bbb Z/2\Bbb Z)$ endowed with the topology of convergence in measure is a Polish space.
Therefore  we will consider  the set $\Cal Z$ of $\Bbb Z/2\Bbb Z$-valued cocycles of $\Cal R_{\theta_1}$ as a  Polish space.
The following  properties of this topological space hold:
\roster
\item"(i)"
the cohomology class of every cocycle is dense in $\Cal Z$,
\item"(ii)"
the subset of ergodic cocycles $\beta$
is a dense $G_\delta$ in $\Cal Z$,
\item"(iii)"
the subset of cocycles $\beta$ such that $(F'_{\theta_1})_\beta$  is rigid
\footnote{We recall that a transformation $S\in\text{Aut}(X,\mu)$ is called rigid if there is a sequence $n_i\to\infty$ such that $S^{n_i}\to\text{Id}$ weakly as $i\to\infty$.} is a $G_\delta$ in $\Cal Z$,
\item"(iv)"
if $\beta(z,F'_{\theta_1}z)=1$ for all $z$ then
$(F'_{\theta_1})_\beta$ is rigid,
\item"(v)"
if $t$ is rationally independent with $\theta_1$ then
 the subset of cocycles $\beta$ such that
$\beta\circ F'_{t}\not\approx\beta$ is residual in $\Cal Z$.
\endroster
The properties (i)--(iii) are well known.
If $\beta(z,F'_{\theta_1}z)=1$ for all $z$ then
$(F'_{\theta_1})_\beta$ has pure point spectrum and hence (iv) follows.
The property (v) follows from \cite{GLS, Theorem~1.2} or \cite{Da1, Theorem~4.2}.

It follows from (i)--(v) that given $t_0>0$ which is rationally independent  with $\theta_1$, there exists a cocycle $\beta:\Cal R_{\theta_1}\to\Bbb Z/2\Bbb Z$ such that the transformation $(F'_{\theta_1})_\beta$ is ergodic, rigid and
$$
L:=\{t\in\Bbb R \mid\beta\circ F_t'\approx\beta\}\not\ni t_0.
$$
Of course, $L$ is a subgroup of $\Bbb R$ and $L\supset\Bbb Z$.
It is well known (for instance, see \cite{New}) that for each transformation $A\in C((F'_{\theta_1})_\beta)$ there are $t\in L$ and a map $\phi\in\Cal M(\Bbb T,\Bbb Z/2\Bbb Z)$ such that
$$
A(z,i)=(F'_tz,i+\phi(z))\quad \text{for all}\quad (z,i)\in\Bbb T\times \Bbb Z/2\Bbb Z.
$$
The map $A\mapsto t+\Bbb Z$ is a group homomorphism from $C((F'_{\theta_1})_\beta)$ onto the quotient group $L/\Bbb Z$.
The kernel of this homomorphism is isomorphic to $\Bbb Z/2\Bbb Z$.
Since $(F'_{\theta_1})_\beta$ is rigid, $C((F'_{\theta_1})_\beta)$ is uncountable.
Hence $L$ is uncountable.

{\bf (B})
Consider the torus $(Y,\nu):=(\Bbb T\times\Bbb T,\lambda_\Bbb T\times\lambda_\Bbb T)$.
Let $Q:=F'_{\theta_1}\times F'_{\theta_2}$,
where  $\theta_2$ is an irrational such that the reals $1, \theta_1, \theta_2$ are rationally independent.
Then $Q$ is a transformation of $(Y,\nu)$ with  pure point spectrum.
Denote by $\Cal R_Q$ the $Q$-orbit equivalence relation.
Then $\Cal R_Q$ is discrete and ergodic.
Define a flow $F=(F_t)_{t\in\Bbb R}$ on $(Y,\nu)$ by setting:
$F_t:=F_t'\times F_{\theta_3t}'$,
where $\theta_3$ is an irrational such that reals $1,\theta_3,\theta_1\theta_3+\theta_2$ are  rationally independent.
It is straightforward  to verify that  $F$ is  strictly $\Cal R_Q$-outer.
Consider now an extension $\beta\otimes 1:\Cal R_Q\to \Bbb Z/2\Bbb Z$ of
$\beta$ given by
$$
\beta\otimes 1(z,w,z',w'):=\beta(z,z').
$$
Then $\beta\otimes 1$ is a cocycle of $\Cal R_Q$.
We first claim that it is ergodic.
Indeed, the $(\beta\otimes 1)$-skew product extension $Q_{\beta\otimes 1}$ of $Q$ is isomorphic to the product $(F_{\theta_1})_\beta\times F_{\theta_2}$ (the corresponding isomorphism is given by a permutation of coordinates in the space $\Bbb T\times\Bbb T\times\Bbb Z/2\Bbb Z$ where
$Q_{\beta\otimes 1}$ acts).
The discrete spectrum of $(F_{\theta_1})_\beta$ is $\{n\theta_1+\Bbb Z\in\Bbb T\mid n\in\Bbb Z\}$.
 Hence it intersects trivially with the  discrete spectrum of $F_{\theta_2}$ which  is $\{n\theta_2+\Bbb Z\in\Bbb T\mid n\in\Bbb Z\}$.
Therefore $(F_{\theta_1})_\beta\times F_{\theta_2}$ is ergodic, as desired.

Next we claim that
$$
\{t\in\Bbb R\mid(\beta\otimes 1)\circ F_t\approx\beta\otimes 1\}=L.\tag2-1
$$
The inclusion $\supset$ is obvious.
To prove the converse, let $f:Y\to\Bbb Z/2\Bbb Z$ be a map such that
$$
\beta\circ F_t'(z,z')=-f(z,w)+\beta(z,z')+f(z',w')
$$
for some $t\in\Bbb R$ and all $(z,w,z',w')\in\Cal R_Q\cap (Y'\times Y')$, where $Y'$ is a $\nu$-conull subset  in $Y$.
Without loss of generality we can think that $Y'$ is $\text{Id}\times F'_{\theta_2}$-invariant.
It follows that
$$
\beta\circ F_t'(z,z')=-f(z,F'_{\theta_2}w)+\beta(z,z')+f(z',F'_{\theta_2}w')
$$
and hence the function
$(z,w)\mapsto f(z,w)-f(z,F'_{\theta_2}w)$ is $\Cal R_Q$-invariant.
Since $\Cal R_Q$ is ergodic, this function is constant.
This implies that $f(z,w)=f(z,F'_{2\theta_2}w)$ for a.a. $(z,w)\in Y$, i.e. $f$ is invariant under the transformation $\text{Id}\times F'_{2\theta_2}$.
Since the transformation $F'_{2\theta_2}$ of $Y$ is ergodic, $f$ does not depend on $w$.
Thus $\beta\circ F'_t\approx\beta$ and hence the inclusion $\subset$ in~\thetag{2-1} is established.

{\bf (C)}
Let $(\widetilde Y,\widetilde\nu):=(Y\times \Bbb T,\nu\times\lambda_\Bbb T)$.
Define an equivalence relation $\Cal T$ on $(\widetilde Y,\widetilde\nu)$ by setting
$$
(y,z)\sim_{\Cal T} (y',z')\quad\text{if}\quad
y\sim_{\Cal R_Q} y'.
$$
Then $\Cal T$ is an amenable ergodic continuous measure preserving equivalence relation on $\widetilde Y$.
Define a flow $\widetilde F=(\widetilde F_t)_{t\in\Bbb R}$ on $\widetilde Y$ by setting
$\widetilde F_t:=F_t\times\text{Id}$, $t\in\Bbb R$.
Then $\widetilde F$ is strictly $\Cal T$-outer.
 Next, consider the cocycle $\widetilde\beta:=\beta\otimes 1\otimes 1:\Cal T\to\Bbb Z/2\Bbb Z$  of $\Cal T$.
Then, of course, $\widetilde\beta$ is ergodic and
$$
\{t\in\Bbb R\mid \widetilde\beta\circ\widetilde F_t\approx\widetilde \beta\}=
\{t\in\Bbb R\mid(\beta\otimes 1)\circ F_t\approx\beta\otimes 1\}=L.
$$
\qed
\enddemo

We will need  an auxiliary fact which is a particular case of  \cite{VF, Theorem~1}.

\proclaim{Lemma 2.3} Let $\Cal R_i$ be an amenable ergodic continuous measure preserving equivalence relation on a standard probability space $(X_i,\goth B_i,\mu_i)$ and let $V^{(i)}=(V^{(i)}_h)_{h\in H}$ be a strictly $\Cal R_i$-outer action of an amenable locally compact second countable group $H$ on $X_i$, $i=1,2$.
Then there is a Borel isomorphism $R:(X_1,\mu_1)\to (X_2,\mu_2)$
and  two conull subsets $Y_1\subset X_1$ and $Y_2\subset X_2$
such that
$$
(R\times R)(\Cal R_1\cap (Y_1\times Y_1))=\Cal R_2\cap (Y_2\times Y_2)
$$
and for each $h\in H$, there exists an $\Cal R_2$-inner transformation $S_h$ of $X_2$ with
$$
RV^{(1)}_hR^{-1}=V^{(2)}_hS_h.
$$
\endproclaim

The following lemma is perhaps well known.
However we were unable to find its proof  in the literature.
Therefore we provide our proof of it.

\proclaim{Lemma 2.4}
Let $H=(H_s)_{s\in \Bbb R}$ and $G=(G_t)_{t\in\Bbb R}$ be the horocycle flow and the geodesic flow on a surface $X$ of constant negative curvature.
 Let $\mu$ denote the normalized volume on $X$.
Then the joint action $\Bbb R\rtimes\Bbb R\ni(s,t)\mapsto H_sG_t$ of
the semidirect product $\Bbb R\rtimes\Bbb R$ on $(X,\mu)$ is free (mod 0).
\endproclaim
\demo{Proof}
We define  multiplication on $\Bbb R\rtimes\Bbb R$ by setting
$$
(s,t)(s',t'):=(s+e^{t}\cdot s',t+t').
$$
Without loss of generality we may assume that $X=\Gamma\backslash\text{SL}_2(\Bbb R)$ for a lattice $\Gamma\subset\text{SL}_2(\Bbb R)$, $\mu$ is  Haar measure on $X$ and
$$
H_s(\Gamma g)=\Gamma g
\left(\matrix 1&0\\s&1\endmatrix\right),\qquad G_t(\Gamma g)=\Gamma g
\left(\matrix e^{t/2}&0\\0& e^{-t/2}\endmatrix\right),
$$
$g\in \text{SL}_2(\Bbb R)$, $t,s\in\Bbb R$  \cite{Ra}.
Then the action $\Bbb R\rtimes\Bbb R\ni(s,t)\mapsto H_sG_t$ is well defined.

Since $H$ is free, we only need to show that
 the subset
$$
\left\{g\in \text{SL}_2(\Bbb R)\bigg|
\, \Gamma g
\left(\matrix
a&0\\b& a^{-1}
\endmatrix\right)=\Gamma g \text{ \ for some }b\in\Bbb R\text{ and } 1\ne a>0\right\}
$$
is of zero Haar measure in $\text{SL}_2(\Bbb R)$.
For this purpose, we will show  that for each $\gamma\in\Gamma$,  the subset
$$
M_\gamma:=\left\{g\in \text{SL}_2(\Bbb R)\bigg|\
g\left(\matrix
a&0\\b&a^{-1}
\endmatrix\right)g^{-1}=\gamma
\text{ \ for some }b\in\Bbb R\text{ and } 1\ne a>0
\right\}
$$
is of zero measure in $\text{SL}_2(\Bbb R)$.
Indeed, given $g_1,g_2\in M_\gamma$, the product $g_1g_2^{-1}$ commutes with the matrix $\left(\matrix
a&0\\b&a^{-1}
\endmatrix\right)$.
Since $a\ne 1$, it follows that $g_1g_2^{-1}$ is a lower-triangular matrix.
It remains to note that  the Haar measure of the subgroup of lower-triangular matrices in $\text{SL}_2(\Bbb R)$ is zero.
\qed
\enddemo

 Let  $T=(T_f)_{f\in F}$ be an action of a locally compact Abelian group  $F$ on $(X,\mu)$.
 A measure $\nu$ on $X\times X$ is a {\it 2-fold self-joining} of $T$ if
 $\nu$ is invariant under the diagonal action $(T_f\times T_f)_{f\in F}$ and the marginal projections of $\nu$ are both equal $\mu$.
If  each ergodic 2-fold self-joining of $T$ is either $\mu\times\mu$ or a measure supported by the graph of $T_f$ for some $f\in F$ then $T$ is said {\it to have property MSJ$_2$} (two-fold minimal self-joinings) \cite{dJR}.

\demo{Proof of Theorem~2.1}
Let $H=(H_s)_{s\in \Bbb R}$ and $G=(G_t)_{t\in\Bbb R}$ be the horocycle flow  and  geodesic flow on the surface $(X,\mu)$ of constant negative curvature.
Suppose that $H$ has the property of MSJ$_2$ (see \cite{Ra}).
 It follows from Lemma~2.4 that there is a Borel $H$- and $G$-invariant $\mu$-conull subset $X_0\subset X$ such that
\roster
\item"$(\bullet)$"
if $G_tH_sx=x$ for some $t,s\in\Bbb R$ and $x\in X_0$ then $t=s=0$.
\endroster
Denote by $\Cal R$ the $H$-orbit equivalence relation on $X_0$.
Then $\Cal R$ is amenable ergodic continuous and $\mu$-preserving.
It follows from $(\bullet)$ that $G$ is strictly $\Cal R$-outer.

Let  $\widetilde Y,\widetilde\nu,\Cal T,\widetilde F,\widetilde\beta,L,t_0$
denote the same objects as in Proposition~2.2.
Then by Lem\-ma~2.3, there is a Borel isomorphism $R:(\widetilde Y,\widetilde\nu)\to(X_0,\mu)$ and  conull subsets $X_1\subset X_0$
and $\widetilde Y_1\subset\widetilde Y$ such that
$(R\times R)(\Cal T\cap(\widetilde Y_1\times \widetilde Y_1))=\Cal R\cap(X_1\times X_1)$
and for each $t\in\Bbb R$,
$$
R\widetilde F_tR^{-1}=G_tS_t,
$$
where $S_t$ is an $\Cal R$-inner transformation of $X_0$.
Denote by $\alpha$ the cocycle $\widetilde\beta\circ R^{-1}:\Cal R\to\Bbb Z/2\Bbb Z$.
Since $\widetilde\beta$ is ergodic, so is  $\alpha$.
 For each $t\in\Bbb R$,
$$
\alpha\circ G_t=\widetilde\beta\circ(\widetilde F_tR^{-1}S_t^{-1}).\tag2-3
$$
It follows  from Proposition~2.2 that the cocycle in the right-hand side of \thetag{2-3} is cohomologous to $\widetilde\beta\circ(R^{-1}S_t^{-1})=\alpha\circ S_t^{-1}$ if and only if $t\in L$.
Since $S_t$ is $\Cal R$-inner, we obtain $\alpha\circ S_t^{-1}\approx\alpha$.
Thus,
$$
\{t\in\Bbb R\mid \alpha\circ G_t\approx\alpha\}=L,
$$
and $L$ is a proper uncountable subgroup of $\Bbb R$.

Denote by $H^\alpha$ the $\alpha$-skew product extension of $H$.
We will show that
$$
I(H^\alpha)\cap\Bbb R_+=\{e^a\mid a\in L\}.\tag2-4
$$
Given  $a\in\Bbb R$, the real  $e^a$ belongs to $I(H^\alpha)$ if and only if
there exists a transformation $V$ of $X\times\Bbb T$ such that $V\bullet H^\alpha=H^\alpha\circ e^a$ (for the definition of $\bullet$ and $\circ$ we refer to Section~1).
Now we are going to describe the ``structure'' of $V$.
Denote by $\kappa$ the corresponding graph-joining of $H^\alpha$ and $H^\alpha\circ e^a$, i.e. $\kappa$ is supported on the graph of $V$.
Thus $\kappa$ is a measure on $X\times\Bbb Z/2\Bbb Z\times X\times\Bbb Z/2\Bbb Z$.
 Denote by $\kappa'$ the projection of $\kappa$ on $X\times X$.
 Then $\kappa'$ is an ergodic joining of $H$ and $H\circ e^a$.
Hence $\kappa'\circ(\text{Id}\times G_a)$ is an ergodic 2-fold self-joining of $H$.
Since $H$ has MSJ$_2$, it follows that either $\kappa'=\mu\times\mu$ or
$\kappa'$ is a graph-joining supported on the graph of $G_aH_s$ for some $s\in\Bbb R$.
In the first case we get a contradiction to the fact that $\kappa$ is a graph-joining.
Therefore the second case holds.
Then
$$
V(x,\cdot)=(G_aH_sx,\cdot).
$$
Replacing $V$ with $VH_{-s}^\alpha$ we can assume without loss of generality that
$$
V(x,\cdot)=(G_ax, \cdot).
$$
It is a standard trick to show that such a $V$ conjugates
$H^\alpha$ with $H^\alpha\circ e^a$ if and only if
$V(x,z)=(G_ax,i+\phi(x))$ for some Borel function $\phi:X\to \Bbb Z/2\Bbb Z$ such that the cocycles $\alpha\circ G_a$ and $\alpha$ are cohomologous, i.e. $a\in L$.
Thus, \thetag{2-4} is established.
Therefore $I(H^\alpha)$ is uncountable and $I(H^\alpha)\ne \Bbb R_+^*$.
As we noted in the beginning of this section, the uncountability of $I(H^\alpha)$ implies that $H^\alpha$ is weakly mixing.
Since $H$ is mixing of all orders, we deduce from \cite{Ru} that   $H^\alpha$ is also mixing of all orders.
\qed
\enddemo

\remark{Remark 2.5} (i) In a similar way, one can obtain the following generalization of
Theorem~2.1.
Let $R$ be an irrational rotation on $(\Bbb T,\lambda_{\Bbb T})$ and let $K$ be a compact second countable group.
Denote by Aut$\,K$ the group of continuous automorphisms of $K$.
Fix an ergodic cocycle $\beta$ of the $R$-orbital equivalence relation $\Cal D$ to $K$.
We let
$$
L(\Cal D,\alpha):=\{S\in C(R)\mid \beta\circ S\approx v\circ\beta, v\in\text{Aut}\,K\}.
$$
Since $C(R)=\Bbb T$, we denote by $\pi:\Bbb R\to C(R)$ the canonical projection $t\mapsto \pi(t):=t+\Bbb Z$.
Then there is a mixing flow $T$ such that
$$
I(T)\cap\Bbb R_+^*=\{e^t\mid \pi(t)\in L(\Cal D,\beta)\}.
$$
Thus we obtain a class of flows $T$ for which the invariant $I(T)$ is of purely ``cohomological'' nature.

(ii) We also note that the groups like $L(\Cal D,\beta)$ and their orbital analogues appear naturally when studying extensions of ergodic dynamical systems and equivalence relations.
For more information about them we refer the reader to  \cite{Da1}, \cite{DaG} and references therein.
\endremark

\head 3. Flows without self-similarity
\endhead

In this section we study  the problem of existence of self-similarities from the Baire category point of view.

We first show that a generic flow has no self-similarities by using some examples of such flows from \cite{FrL}.
Then an alternative, independent from \cite{FrL}, proof of this fact is given.
We construct explicitly a  {\it rank-one} flow such that
the dilations of the measure of maximal spectral type of this flow are mutually orthogonal.
This property is generic in Flow$(X,\mu)$.
It follows from the existence of certain special {\it weak limits} of the flow (see \thetag{3-1} and \thetag{3-2} below).
To manufacture these weak limits we combine two standard  rank-one constructions with {\it flat} and {\it staircase} roofs.
As a corollary, we obtain that
a generic transformation does not embed into a flow with self-similarities.

Let $\Cal P$ stand for the set of continuous probability measures on the one-point compactification $\overline{\Bbb R}=\Bbb R\cup\{\infty\}$ of $\Bbb R$.
Then $\Cal P$ is a compact metric space in the $*$-weak topology.
Denote by $\Cal P_C\subset\Cal P$ the subset of continuous, i.e. non-atomic, measures.
It is well known that $\Cal P_C$ is a dense $G_\delta$ in $\Cal P$ \cite{Na}.
Hence it is Polish when endowed with the induced topology.
Since $\sigma(\{\infty\})=0$ for each $\sigma\in\Cal P_C$,
we  identify $\Cal P_C$ with the space of non-atomic probability measures on $\Bbb R$.

Given $\sigma\in \Cal P$ and $t\ne 0$, we define a measure $\sigma_t$ by setting
$\sigma_t(A):=\sigma(t\cdot A)$ for each Borel subset $A\subset \Bbb R$.

\proclaim{Lemma 3.1}
The set
$\Cal S:=\{\sigma\in \Cal P_C\mid \sigma_t\perp\sigma \text{ for all }t>0, t\ne 1\}$
is a $G_\delta$ in $\Cal P_C$.
\endproclaim
\demo{Proof}
For each open subset $O\subset \Bbb R$ with the finite boundary, the map
$$
\Sigma:\Cal P_C\times\Bbb R_+^*\ni (\sigma,t)\mapsto\sigma_t(O)\in\Bbb R
$$
is continuous.
Therefore, given a segment $I\not\ni 1$ in $\Bbb R_+^*$ and an open subset $O\subset \Bbb R$, the map
$$
f_{O,I}: \Cal P_C\ni\sigma\mapsto (\sigma(O),\max_{t\in I}\sigma_t(O))\in\Bbb R^2
$$
is continuous.
Given a segment $I\subset\Bbb R$ and $N>0$, we denote by $\Cal P_N(I)$ the partition of a segment $I$  into $N$ sub-segments of equal length.
We let
$$
\Cal S':=\bigcap_{I\not\ni 1}\bigcap_{n\in\Bbb N}\bigcup_{N>0}\bigcap_{\Delta\in\Cal P_N(I)}\bigcup_{O}f_{O,\Delta}^{-1}((1-1/n,+\infty)\times(-\infty,1/n)),
$$
where $I$ runs over  segments with positive rational endpoints and $O$ runs over the collection of  open subsets in $\Bbb R$ with the finite boundary.
Of course, $\Cal S'$ is $G_\delta$ in $\Cal P_C$.
It is easy to see that $\Cal S'\subset\Cal S$.
Let us show the converse inclusion.
Indeed, if $\sigma\in\Cal S$ then for each  $t\in\Bbb R^*_+$ and $n\in\Bbb N$, there is an open subset $O_{n,t}\subset\Bbb R$ with the finite boundary such that
$\sigma(O_{n,t})>1-1/n$ and $\sigma_t(O_{n,t})<1/n$.
Since $\Sigma$ is continuous,
for each $t\in\Bbb R^*_+$ there is a neighborhood $U(t)$ of $t$ such that
$\sigma_\tau(O_{n,t})<1/n$ for all $\tau\in U(t)$.
Therefore, we can assume without loss of generality that for each segment $I\not\ni 1$, there exists $N>0$ such that the map $t\mapsto O_{n,t}$ is constant for every subsegment $\Delta\in\Cal P_N(I)$.
Hence,
$$
\sigma\in\bigcap_{\Delta\in\Cal P_N(I)}\bigcup_{O}f_{O,\Delta}^{-1}((1-1/n,+\infty)\times(-\infty,1/n)).
$$
Thus $\Cal S=\Cal S'$.
\qed
\enddemo

The subset $\Cal W$ of weakly mixing flows on $(X,\goth B,\mu)$ is a dense $G_\delta$ in $\text{Flow}(X,\mu)$ because of the following three facts:
\roster
\item"---"
the mapping Flow$(X,\mu)\ni T\mapsto T_1\in\text{Aut}(X,\mu)$ is continuous,
\item"---" the subset of weakly mixing transformations is a $G_\delta$ in
$\text{Aut}(X,\mu)$,
\item"---"
the Aut$(X,\mu)$-orbit\footnote{We mean the action defined by the formula \thetag{1-2}.}  of each weakly mixing flow $T$  is dense in Flow$(X,\mu)$. (As in the case of $\Bbb Z$-actions, this fact follows easily from the Rokhlin lemma.
 See e.g. \cite{OW} and \cite{DaSo} for more general versions of Rokhlin lemma.)
\endroster
Fix an orthonormal basis $(v_j)_{j\in\Bbb N}$ in $L^2_0(X,\mu)$.
Given $T\in\text{Flow}(X,\mu)$, let $U_T=(U_T(t))_{t\in\Bbb R}$
denote the corresponding {\it Koopman} unitary representation of $\Bbb R$ in $L^2_0(X,\mu)$.
Recall that $U_T(t)f:=f\circ T_{-t}$.
For each $j$, let $\sigma_{T,j}$ be the only probability measure on $\Bbb R$ such that for each $t\in\Bbb R$,
$$
\langle U_T(t)v_j,v_j\rangle=\int_{\Bbb R}\exp(2\pi i\lambda t)\,d\sigma_{T,j}(\lambda).
$$
We now let $\sigma_T:=\sum_{j=1}^\infty 2^{-j}\sigma_{T,j}$.
Then $\sigma_T$ is a measure of the maximal spectral type of $U_T$
and the map
$$
\Cal W\ni T\mapsto\sigma_T\in\Cal P_C
$$
is continuous.

\proclaim{Theorem 3.2}
The subset $\Cal T:=\{T\in\Cal W\mid \sigma_T\in\Cal S\}$ is a dense $G_\delta$ in $\text{\rm Flow}(X,\mu)$.
\endproclaim
\demo{Proof}
 It follows from Lemma 3.1 that $\Cal T$ is $G_\delta$ in $\Cal W$.
In \cite{FrL}, a Gaussian flow $T$ with a simple spectrum was constructed such that $T\in\Cal T$ (see also another example in Proposition~3.4 below).
It remains to use the fact that the  Aut$(X,\mu)$-orbit of each ergodic flow in $\Cal W$ is dense in $\text{\rm Flow}(X,\mu)$.
\qed
\enddemo

We recall two concepts of disjointness for dynamical systems.
 Let we are given two actions $T=(T_a)_{a\in A}$ and $S=(S_a)_{a\in A}$ of a locally compact  second countable  Abelian group $A$ on standard probability spaces $(X,\mu)$ and $(Y,\nu)$ respectively.
 The actions  are called
\roster
\item"(i)" {\it disjoint in the sense of Furstenberg} if $\mu\times\nu$ is the only $(T_a\times S_a)_{a\in A}$-invariant measure on $X\times Y$ with marginals $\mu$ and $\nu$.
We will denote this by $T\perp^{\text F} S$.
\item"(ii)" {\it spectrally disjoint}
if the  maximal spectral types of $T$ and $S$ are mutually orthogonal.
\endroster
If $T$ and $S$ are spectrally disjoint then  $T\perp^{\text F} S$.
The converse is not true.
We note also that if $t=-1$ then $T$ and $T\circ t$ have the same  maximal spectral type.
If $T$ and $S$ are weakly mixing and $A_0$ is a cocompact subgroup in $A$ then $T\perp^{\text F} S$ if and only if $(T\restriction A_0)\perp^{\text F} (S\restriction A_0)$ \cite{dJR}.

\proclaim{Corollary 3.3}   For a generic flow $T$, the group $I(T)$ is trivial.
Moreover, $\sigma_T\perp\sigma_{T\circ t}$ if  $|t|\ne 1$ and
$T\perp^{\text F} T\circ t$ if $t=-1$.
\endproclaim
\demo{Proof}
(i) Take $T\in\Cal T$.
Since $\sigma_{T\circ t}=(\sigma_T)_{t}$ and $(\sigma_T)_{-t}\sim (\sigma_T)_{t}$,  we obtain
that
the flows $T\circ t$ and $T$ are spectrally disjoint for all and $t\in\Bbb R^*$, $t\ne -1$.
Hence $I(T)\subset\{-1,1\}$.

(ii)The map $\text{Flow}(X,\mu)\ni T\mapsto T_1\in\text{Aut}(X,\mu)$
is continuous.
By \cite{dJ}, the set of transformations $S$  such that $S\perp^{\text{F}}S^{-1}$ is a dense $G_\delta$ in $\text{Aut}(X,\mu)$.
An example a weakly mixing flow $T$ with $T_1\perp^{\text{F}}T_{-1}$  was given in \cite{dJP}.
It follows that the set
$\Cal A:=\{T\in\Cal W\mid T_1\perp^{\text{F}}T_{-1}\}$  is a dense $G_\delta$ in $\Cal W$.
If $T\in\Cal A$ then  $-1\not\in I(T)$.
\qed
\enddemo

We now give an explicit example of a rank-one flow $T\in\Cal T$.
For that we  recall a classical cutting-and-stacking construction of  rank-one flows.
The construction process is inductive.
Suppose we are given
\roster
\item"(a)" a sequence of integers  $r_n>1$, and
\item"(b)" a sequence of mappings
$s_n:\{1,\dots,r_n\}\to\Bbb R_+$.
\endroster
On the $n$-th step  we have a tower, say $X_n$, which is a rectangular of height $h_n$ and width $w_n$.
We cut it into $r_n$ subtowers of equal width $w_n/r_n$.
 Enumerate these subtowers  from the left to the right by $1,\dots,r_n$.
Then for each $j=1,\dots,r_n$, we put a rectangle of height $s_n(j)$ and
width $w_{n+1}:=w_n/r_n$ on the top of the $j$-th subcolumn.
Thus we obtain a family of $r_n$ enumerated towers of height
$$
h_n+s_n(1),h_n+s_n(2),\dots,h_n+s_n(r_n).
$$
All of them have the same width $w_{n+1}$.
We now stack these towers in the following way: put the  second tower on the top of the first tower, the third tower on the top of the second one and so on.
Then we obtain a new tower $X_{n+1}$ of height $h_{n+1}:=r_nh_n+\sum_{j=1}^{r_n}s_n(j)$ and width $w_{n+1}$.
Since $X_{n+1}$ is embedded into $\Bbb R^2$, we endow it with the induced Lebesgue measure, say $\mu_{n+1}$.

Continuing this procedure infinitely many times we obtain  a $\sigma$-finite standard non-atomic measure space $(X,\mu)$ as an inductive limit of the sequence of finite measure spaces $(X_0,\mu_0)\subset(X_1,\mu_1)\subset\cdots$.
It is easy to see that $\mu$ is finite if and only if
$$
\sum_{n=1}^\infty h_n^{-1}r_n^{-1}\sum_{j=1}^{r_n}s_n(j)<\infty.
$$
We will say that a function $f: X\to \Bbb C$ is $X_n$-{\it measurable}
if $f$ is supported on $X_n$ and $f(x)=f(x')$ whenever $x$ and $x'$ are on the same height in $X_n$.

We now define a flow $T=(T_t)_{t\in\Bbb R}$ on $X$ by setting
$$
T_t(y,z):=(y, t+z), \quad\text{whenever } (y,z), (y,t+z)\in X_{n},
$$
$n=0,1,\dots$.
Geometrically this means that $T_t$ moves a point in $X_n$ up with a unit speed until the point reaches the  top of $X_n$.
It is easy to verify that $T$ is well defined  on the entire space (more precisely, on a $\mu$-conull subset of) $X$  when $n\to\infty$.
This flow preserves $\mu$.
We call $T$ the {\it rank-one flow} associated with $(r_n,s_n)_{n=1}^\infty$.

\proclaim{Proposition 3.4}
Let $T$ be a finite measure preserving rank-one flow associated with a sequence $(r_n,s_n)_{n=1}^\infty$ and let $r_n=10^{n}$ for all $n$.
Suppose that there are a sequence of positive integers $n_k\to\infty$
and a sequence of positive reals $u_k\to 0$ such that
$u_kr_{n_k}\to\infty$ and for each $k$,
\roster
\item"\rom{(i)}"
$s_{n_k-1}(j)=0$ for all $1\le j\le r_{n_k-1}$ and
\item"\rom{(ii)}"
$s_{n_k}(j)=(j-1)u_k$ for all $1\le j\le r_{n_k}$.
\endroster
Then $T\in\Cal T$.
\endproclaim

\demo{Proof}
The conditions (i) and (ii) mean that  infinitely many towers, numbered with $n_k-1$, have  a {\it flat} roof with no spacers added at all while the subsequent towers, numbered with $n_k$, have a  {\it staircase} roof.

Fix $l>0$.
We claim that
$$
\align
&U_T(-dh_{n_k})\to 10^{-l}I\quad \text{for $d=1-10^{-l}$ and}\tag 3-1\\
&U_T(-ch_{n_k})\to 0\quad \text{ uniformly in } c\in [1,10^l],\tag3-2
\endalign
$$
where the arrows mean the convergence in the weak operator topology as $k\to\infty$.
It follows from \thetag{3-1} and \thetag{3-2} and the spectral theorem for $U_T$ that
 $\sigma_{T\circ d}\perp\sigma_{T\circ c}$ for all
$c\in [1,10^l]$.
Hence $\sigma_T\perp(\sigma_T)_t$ for all $t\in(d^{-1},10^l)$.
Since $l$ is arbitrary, we obtain
$\sigma_T\in\Cal S$.
This implies easily that $T\in\Cal W$.
Hence $T\in\Cal T$.

It remains to prove \thetag{3-1} and \thetag{3-2}.
Another piece of notation will be needed.
Given a function $f\in L^2(X,\mu)$,  denote by $f_{k,i}$ the restriction of $f$ to the $i$-th subtower of $X_{k}$, $1\le i\le r_{k}$, i.e.
$f_{k,i}(x)=f(x)$ if $x$ belongs to the $i$-th subtower and $f_{k,i}(x)=0$ otherwise.

First of all we verify $U_T(-h_{n_k})\to 0$.
Take $f$ in the unit ball of $L^2_0(X,\mu)$.
Then for each $\epsilon>0$, there is $k_0>0$ and  $f'\in L^2_0(X,\mu)$ such that $f'$ is $X_{k_0}$-measurable, $\|f-f'\|_2\le \epsilon$,  and $|f'|<D$ for some real $D$.
Take $k$ such that $n_k>k_0$.
Cross out from $X_{n_k}$  the bottom layer of height $(r_{n_k}-2)u_k$.
Denote the rest of $X_{n_k}$ by $X_{n_k}^0$.
Since $f'$ is bounded and $\mu(X_{n_k})-\mu(X_{n_k}^0)\to 0$,
we can assume without loss of generality that $f'$ is supported on $X_{n_k}^0$.
 Then $f'=\sum_{j=1}^{r_{n_k}}f'_{n_k,j}$.
It is easy to deduce from (ii) that $f'_{n_k,j}\circ T_{h_{n_k}}=f'_{n_k,j+1}\circ T_{-u_kj}$ for all $1\le j< r_{n_k}$.
We have
$$
\align
\langle U_T(-h_{n_k})f,f\rangle &=\sum_{j,q=1}^{r_{n_k-1}}\langle f_{n_k,j+1}'\circ T_{-ju_k},f'_{n_k,q}\rangle\pm 2\epsilon\pm2\|f'_{n_k,r_{n_k}}\|_2\\
&=\sum_{j=1}^{r_{n_k-1}}\langle (f'\circ T_{-ju_k})_{n_k,j+1},f'_{n_k,j+1}\rangle\pm 2\epsilon\pm\frac{2D}{r_{n_k}}\\
&=\sum_{j=1}^{r_{n_k-1}}\frac 1{r_{n_k}}\langle U_T(ju_k)f',f'\rangle\pm 2\epsilon\pm\frac{2D}{r_{n_k}}\\
&=\bigg\langle\bigg(\frac 1{r_{n_k}}\sum_{j=1}^{r_{n_k-1}} U_T(ju_k)\bigg)f,f\bigg\rangle\pm 4\epsilon\pm\frac{2D}{r_{n_k}}.
\endalign
$$
Applying the mean ergodic theorem  we obtain that
$\langle U_T(-h_{n_k})f,f\rangle\to 0$, as desired.
Only a slight modification of the above argument is needed to prove
the following fact: for each integer $p>0$,
$$
\sup_g|\langle U_T(-ph_{n_k})f,g\rangle|\to 0,\tag3-3
$$
where the supremum is taken over all $X_{n_k}$-measurable functions $g$ with $\|g\|_2\le 1$.

Now let $f'$ be an $X_{n_k-1}$-measurable bounded function.
Since $h_{n_k}=h_{n_k-1}r_{n_k-1}$,
it follows from (i) that
$$
f'_{n_k-1,j}\circ T_{dh_{n_k}}=f'_{n_k-1,j+dr_{n_k-1}}\quad\text{for all }1\le j\le 10^{-l}r_{n_k-1}.
$$
We let  $f^\bullet:=\sum_{j=1}^{10^{-l}r_{n_k-1}}f'_{n_k-1,j}$
and $f^\circ:=f'-f^\bullet-f'_{n_k-1,r_{n_k-1}}$.
Then
$$
\aligned
\langle U_T(-dh_{n_k})f^\bullet,f'\rangle &=\sum_{j=1}^{10^{-l}r_{n_k-1}}
\langle f'_{n_k-1,j+dr_{n_k-1}},f'\rangle\\
&=\sum_{j=1}^{10^{-l}r_{n_k-1}}\frac {\|f'\|^2_2}{r_{n_k-1}}=10^{-l}\|f'\|^2_2
\endaligned
\tag3-4
$$
and
$$
\langle
U_T(-dh_{n_k})f^\circ,f'\rangle
=
\langle f^\circ\circ T_{(d-1)h_{n_k}},U_T(h_{n_k})f'
\rangle.
\tag3-5
$$
It is easy to verify that the function $f^\circ\circ T_{(d-1)h_{n_k}}$ is $X_{n_k}$-measurable.
Therefore it follows from \thetag{3-3} and \thetag{3-5} that $\langle
U_T(-dh_{n_k})f^\circ,f'\rangle\to 0$.
This fact plus \thetag{3-4} imply \thetag{3-1}.

To show \thetag{3-2} we take $c\in[1,10^l]$ and write $ch_{n_k}$ as $ch_{n_k}=c_kh_{n_k}+c_k'$ with $c_k\in\Bbb N$ and $0\le c_k'< h_{n_k}$.
Partition $X_{n_k}$ by a horizontal line on the height $h_{n_k}-c_k'$ into two subsets $X_{n_k}^0$ (bottom part) and $X_{n_k}^1$ (upper part).
Take a bounded $X_{n_k}$-measurable function $f'$.
Then
$$
\align
&\langle U_T(-ch_{n_k})f',f' \rangle\\
&=
 \langle U_T(-c_kh_{n_k})(f'1_{X^0_{n_k}})\circ T_{c_k'},f' \rangle
+
\langle U_T(-(c_k+1)h_{n_k})(f'1_{X^1_{n_k}})\circ T_{c_k'-h_{n_k}},f' \rangle\\
&=\langle (f'1_{X^0_{n_k}})\circ T_{c_k'},U_T(c_kh_{n_k})f' \rangle
+
\langle (f'1_{X^1_{n_k}})\circ T_{c_k'-h_{n_k}},U_T((c_k+1)h_{n_k})f' \rangle.
\endalign
$$
Since the functions $f'1_{X^0_{n_k}}\circ T_{c_k'}$ and
$(f'1_{X^1_{n_k}})\circ T_{c_k'-h_{n_k}}$ are $X_{n_k}$-measurable, we
can apply \thetag{3-3} to obtain $\sup_{1\le c\le 10^{l}}|\langle U_T(-ch_{n_k})f',f' \rangle|\to 0$ as $k\to\infty$.
\qed
\enddemo

\remark{Remark \rom{3.5}}
We note that it follows directly from Proposition 3.4 that $\Cal T$ is residual.
Indeed, the subset $\Cal L$ of flows $T$ such that for each $l>0$, the
limits
\thetag{3-1} and \thetag{3-2} exist along a common subsequence of $(n_k)_{k>1}$ is a $G_\delta$ in Flow$(X,\mu)$.
This subset is invariant under the action of Aut$(X,\mu)$ by  conjugation. Hence if it is non-empty then it is dense.
As  was shown in the proof of Proposition~3.4, $\emptyset\ne\Cal L\subset\Cal T$.
Thus $\Cal T$ is residual.
On the other hand, the statement of Theorem~3.2 (which uses Lemma~3.1) is sharper: $\Cal T$ is $G_\delta$ itself.
\endremark

\proclaim{Theorem 3.6}
\roster
\item"(i)"
 A generic transformation from $\text{\rom{Aut}}(X,\mu)$ embeds into a flow $T$ such that  $I(T)=\{1\}$.
Moreover, it embeds into a flow possessing all the properties listed in Corollary~3.3.
\item"(ii)"
 A generic transformation from $\text{\rom{Aut}}(X,\mu)$ does not embed into a flow $T$ with $I(T)\not=\{1\}$.
\endroster
\endproclaim
\demo{Proof}
(i) It was shown in \cite{dRdS} that the image of a  non-meager subset in Flow$(X,\mu)$ under the map $T\mapsto T_1$ is non-meager in Aut$(X,\mu)$.
If a non-meager subset of Aut$(X,\mu)$  is invariant under the conjugacy then it is residual in Aut$(X,\mu)$ \cite{GK}.
In view of that, (i) follows from Theorem~3.2 and Corollary~3.3.

 (ii) A similar reasoning yields that the set $\{T_1\mid T\in\Cal L\}$ is residual in Aut$(X,\mu)$.
See Remark~3.5 for the definition of $\Cal L$.
 Hence the intersection
$$
\Cal J:=\{T_1\mid T\in\Cal L\}\cap\{S\in\text{Aut}(X,\mu)\mid S\perp^{\text{F}} S^{-1}\text{ and $S$ has a simple spectrum}\}
$$
is also residual in Aut$(X,\mu)$.
Take $J\in\Cal J$ and suppose that $J=Q_1$ for a flow
$Q\in\text{Flow}(X,\mu)$.
 Since $J=T_1$ for  a  flow $T\in\Cal L$ and
  $J$ has a simple spectrum, the flows $T$ and $Q$ commute.
Hence the flow $P:\Bbb R\ni t\mapsto P_t:=T_tQ_t^{-1}$ is well defined.
This flow is periodic, i.e.  $P_{t+1}=P_t$ for all $t\in\Bbb R$.
Since $T\in\Cal L$, we have that for each $l>0$, \thetag{3-1} and \thetag{3-2} hold along a common subsequence of $(n_k)_{k=1}^\infty$.
Therefore  utilizing the fact that
the group $\{P_t\mid t\in\Bbb R\}$ is compact we can pass to a further subsequence,  say $(m_{l,k})_{k=1}^\infty$,   such that
\roster
\item"(a)"
$U_Q(-dh_{m_{l,k}})\to 10^{-l}U_P(\xi)$ for some $0\le\xi<1$ and
for $d=1-10^{-l}$ and
\item"(b)"
$U_Q(-ch_{m_{l,k}})\to 0$ uniformly in  $c\in [1,10^l]$.
\endroster
The conditions (a) and (b) imply that $Q\in\Cal T$ in the same way as
\thetag{3-1} and \thetag{3-2} imply $T\in\Cal T$ in the proof of Proposition~3.4.
Hence $I(Q)\subset\{-1,1\}$.
It follows from the definition of $\Cal J$ that  $Q_1\perp Q_{-1}$.
Therefore
 $-1\not\in I(Q)$.
Thus $I(Q)=\{1\}$.
\qed
\enddemo

\head
 4. Countable groups of self-similarities
\endhead

In this section we construct flows $T$ with a prescribed  countable group $I(T)$.
We solve (Q6), remove a redundant condition from \cite{FrL, Theorem~9.4}
and provide examples of asymmetric (as well as symmetric) Poisson flows.

\subhead I. Rank one and self-similarities
\endsubhead
Let $S$ be a countable subgroup of $\Bbb R_+^*$ such that $S$ considered as a subset
of the $\Bbb Q$-linear space $\Bbb R$ is independent.
It was shown in \cite{FrL} that there is a Gaussian flow $T$ with a simple spectrum
such that $I(T)=S\sqcup(-S)$ and $\sigma_T\perp(\sigma_T)_t$ for each $t\not\in S\sqcup(-S)$.
We prove the existence of a rank-one flow with similar (but not identical) properties.

\proclaim{Theorem 4.1}Let $S$ be a countable subgroup of $\Bbb R_+^*$ such that $S$ considered as a subset
of the $\Bbb Q$-linear space $\Bbb R$ is independent.
Let $T^S$ denote the Cartesian product flow $\bigotimes_{s\in S}T\circ~s$ acting on the space $(X,\mu)^S$.
For a generic flow $T\in\text{\rom{Flow}}(X,\mu)$,%%%%%%%%%%%%%%%%%%%%%%%%%%%%%%%%%%%%%%%%%%%%%%%%%%%%%%%%%%%%%%%
\roster
\item"\rom(i)" the flow
$T^S$ is rank one rigid and weakly mixing,
\item"\rom(ii)"
 $I(T^S)=S$ and, moreover,
\item"\rom(iii)"
$T^S\perp^{\text{\rom F}} (T^S)\circ t$ for each real $t\not\in S\cup\{0\}$.
\endroster
\endproclaim
\demo{Proof}
Since the set of transformations of rank one is a $G_\delta$ in Aut$(X,\mu)$, it follows that the following sets
$$
\align
\Cal O &:=\{T\in\text{Flow}(X,\mu)\mid T_1\text{ is rank one and rigid}\},\\
\Cal O_S &:=\{T\in\text{Flow}(X,\mu)\mid T^S\in O\}
\endalign
$$
are $G_\delta$ in $\text{Flow}(X,\mu)$.
The two sets are non-empty  because they contain any flow with pure point rational spectrum.
Hence they are dense
  in Flow$(X,\mu)$.
Take a flow $T\in \Cal O_S\cap \Cal T$.
Since $T\in\Cal O_S$, the transformation $(T^S)_1$ is  rank one
and rigid.
Hence  $T^S$ is also  rank one and rigid.
It is obvious that $I(T^S)\supset S$.
Now take $r\notin S\cup\{0\}$.
Suppose that  $T^S\not\perp^{\text{F}}T^S\circ r$.
 Then there is an ergodic joining $\rho\ne\mu^S\times\mu^{S}$ of $T^S$ and $T^S\circ r$.
  In other words, $\rho$ is  a measure on $X^{S\sqcup rS}$ which is  invariant under  $\bigotimes_{s\in S\sqcup rS}T\circ s$, and
 the projections of $\rho$ on $X^S$ and $X^{rS}$ are both $\mu^S$.
Since the spectral disjointness implies the disjointness in the sense of Furstenberg, Corollary~3.3 yields that $\rho$ is pairwise independent, i.e. the projection of $\rho$ on any ``coordinate plane'' $X\times X$ is $\mu\times\mu$.
The measure $\sigma_{T\circ z}$ is singular to Lebesgue measure for  each $z\in\Bbb R^*$.
Hence the maximal spectral type of the transformation $(T\circ z)_1=T_z$
is also singular.
Therefore we may apply  Host theorem \cite{Ho}
to the dynamical system $(X^{S\sqcup rS},\rho,\bigotimes_{s\in S\sqcup rS}T_s)$.
 It yields  $\rho=\mu^{S\sqcup rS}$, a contradiction.
Thus  the claims (i)---(iii) are all proved.
\qed
\enddemo

\remark{Remark 4.2}
We note that the condition on $S$ can not be removed from the statement of Theorem~4.1.
 The theorem does not hold whenever $S$ contains a pair of rationally dependent reals.
 This follows from the fact that if $n$ is a positive integer and $T$ is an ergodic flow then the product flow
$T\times T\circ n$ is never of rank one.\footnote{An analogous assertion for $\Bbb Z$-actions was proved by the second named author in \cite{Ry7}.}
 We will show more: {\it the weak closure theorem} (see \cite{Ry2})  does not hold for  this flow,
i.e. the centralizer $C(T\times T\circ n)$ of this flow is not the weak closure of the group $\{T_t\times T_{nt}\mid t\in\Bbb R\}$ in Aut$(X\times X,\mu\times\mu)$.
Indeed, suppose that the weak closure theorem holds for $T\times T\circ n$.
Fix $t>0$.
Since the transformation $\text{Id}\times T_t$ commutes with  $T\times T\circ n$, it follows  that there is a sequence $t_i\to\infty$ such that $T_{t_i}\times T_{nt_i}\to \text{Id}\times T_t$ as $i\to\infty$.
Then on the one hand $T_{t_i}\to\text{Id}$ and hence
 $T_{nt_i}=T_{t_i}^n\to\text{Id}$ but on the other hand
$T_{nt_i}\to T_t\ne\text{Id}$, a contradiction.
\endremark

For  arbitrary countable subgroups $S\subset\Bbb R^*$, we have been unable to find a rank-one weakly mixing flow $T$ with $I(T)=S$.
However we can prove the following (weaker) assertion.

\proclaim{Theorem 4.3}
Let $S$ be a countable subgroup of $\Bbb R^*$.
There is a weakly mixing rank-one rigid flow $T$ such that $I(T)\supset S$.
\endproclaim
\demo{Proof}
Let $\Bbb R\rtimes S$  denote the semidirect product $\Bbb R$ with $S$ with the multiplication as follows:
$$
(r,s)(r',s'):=(r+s\cdot r', ss').
$$
We furnish $\Bbb R\rtimes S$ with  the natural (product) locally compact second countable topology.
Let $\Cal A$ stand for the set of all measure preserving actions of $\Bbb R\rtimes S$ on $(X,\goth B,\mu)$.
We endow $\Cal A$ with the topology of
 uniform convergence on the compacts in $\Bbb R\rtimes S$.
Then $\Cal A$ is a Polish space.
One can show in a standard way that each of the following subsets is a $G_\delta$ in $\Cal A$ :
\roster
\item"(a)"
$\Cal A_1:=\{W\in\Cal A\mid \text{the action $\Bbb R\ni t\mapsto W_{(t,0)}$ is weakly mixing}\}$,
\item"(b)" $\Cal A_2:=\{W\in\Cal A\mid \text{the transformation $ W_{(1,0)}$ is rigid and of rank one}\}$,
\endroster
The two sets are invariant under the action of $\text{Aut}(X,\mu)$ on $\Cal A$ by conjugacy.
Again, using the Rokhlin lemma for $(\Bbb R\rtimes S)$-actions one can show that the conjugacy class of each free $(\Bbb R\rtimes S)$-action
is dense in $\Cal A$.
Therefore $\Cal A_1$ and $\Cal A_2$ are dense $G_\delta$ if they contain at least one free $(\Bbb R\rtimes S)$-action.
Of course, $\Cal A_1$ contains such an action.

It remains to construct a free $(\Bbb R\rtimes S)$-action belonging to  $\Cal A_2$.
Let $\Gamma$ be a dense countable subgroup in $\Bbb R$ such that $S\cdot\Gamma=\Gamma$.
Let $T$ be an ergodic  flow with pure point spectrum $\Gamma$.
Denote by $\widehat\Gamma$ the   Abelian group dual to $\Gamma$.
Then $\widehat\Gamma$ is compact and connected.
We note that  $T$ is defined  on $(\widehat\Gamma,\lambda_{\widehat\Gamma})$ in the following way:
$$
T_t x:=x+h(t),\tag4-1
$$
where $h:\Bbb R\to \widehat\Gamma$ is a continuous  one-to-one homomorphism with dense range in $\widehat\Gamma$.
Of course, $S$ also acts on $\widehat\Gamma$ as follows
$$
s\cdot x(\gamma):=x(s^{-1}\cdot\gamma), \qquad \gamma\in\Gamma.\tag4-2
$$
The two actions \thetag{4-1} and \thetag{4-2} generate a measure preserving  action, say $W$, of $\Bbb R\rtimes S$.
It is easy to verify that $W$ is free and $W\in\Cal A_2$.
\qed
\enddemo

 We note that Theorem~4.3 refines \cite{Ag} and \cite{Da2, Theorem~1.3}.

\subhead II. Poisson suspensions  and Gaussian flows with countable set of self-similarities
\endsubhead
Let $T$ be a measure preserving flow on an infinite $\sigma$-finite measure space $(X,\goth B,\mu)$.
Since the non-trivial constant functions are not integrable,
 the associated Koopman representation $U_T$ is defined on the entire space $L^2(X,\mu)$.
We will always assume that  $T$ has no non-trivial invariant subsets of finite positive measure.
Then the maximal spectral type $\sigma_T$
of $T$ is continuous.
For $t\in\Bbb R$, we denote by $\widetilde T_t$ the Poisson suspension of $T_t$ (see  \cite{Ne} and \cite{Ro}).
Then $\widetilde T:=(\widetilde T_t)_{t\in\Bbb R}$ is a weakly mixing finite measure preserving flow.
 As we noted in \cite{DaR}, if $T$ has a simple spectrum then the {\it Gaussian} flow associated with $\sigma_T$ is spectrally equivalent to $\widetilde T$, i.e. the Koopman representations generated by the two flows are unitarily equivalent.

\proclaim{Theorem 4.4} Let $S$ be a countable subgroup of $\Bbb R^*_+$.
There is a weakly mixing Poisson   suspension flow $\widetilde W$ with a simple spectrum such that $I(\widetilde W)\cap \Bbb R^*_+=S$ and $\sigma_{\widetilde W}\perp(\sigma_{\widetilde W})_t$ for each positive $t\not\in S$.
Hence  there is also a weakly  mixing Gaussian flow $F$
with a simple spectrum such that
$I(F)=S\sqcup(-S)$ and $\sigma_F\perp(\sigma_F)_t$ for each  $t\not\in S\sqcup(-S)$.
\endproclaim

Given a unitary operator $V$ in a  Hilbert space $\Cal H$, we denote by WCP$(V)$ the {\it weak closure of the powers} of $V$, i.e. the closure of the group $\{V^n\mid n\in\Bbb Z\}$ in the weak operator topology.
The unitary operator $\bigoplus_{n\ge 0}V^{\odot n}$ acting in a Hilbert space $\bigoplus_{n\ge 0}\Cal H^{\odot n}$ is called the {\it exponent} of $V$.
It is denoted by $\exp(V)$.

The following two lemmata are well known.
For their proof we refer the reader to, e.g., \cite{Ry6} and \cite{DaR} respectively.

\proclaim{Lemma 4.5} Let $V$ has a simple spectrum.
If
$$
\text{\rom{WCP}}(V)\supset
\{\alpha_nI+\beta_nV\mid n\in\Bbb N\text{ and }\alpha_i/\beta_i\ne\alpha_j/\beta_j\text{ whenever }i\ne j\}
$$
then $\exp(V)$ has a simple spectrum.
\endproclaim

\proclaim{Lemma 4.6}
Let $U, V$ be two unitary operators in a  Hilbert space $\Cal H$.
If $U$ and $V$ have  a  simple spectrum and
 \,$\text{\rom{WCP}}(U\otimes V)\ni aI\otimes V$
  for some $a>0$ then  the tensor product $U\otimes V$ has a simple spectrum.
\endproclaim

We also need the following lemma.

\proclaim{Lemma 4.7}
Let $U=(U(t))_{t\in\Bbb R}$ be a weakly continuous unitary representation of $\Bbb R$ in a  Hilbert space $\Cal H$.
If $U$ has a simple spectrum and \rom{WCP}$(U(c))\ni U(a)$ for some $c,a>0$ with $c/a\not\in\Bbb Q$ then the operator $U(c)$ has a simple spectrum.
\endproclaim
\demo{Proof}
Let $h\in\Cal H$ be a cyclic vector for $U$.
Denote by $\Cal Z$ the $U(c)$-cyclic space generated by $h$.
Since $\Cal Z$ is invariant under any operator from
\rom{WCP}$(U(c))$, it follows that $U(a)h\in\Cal Z$.
By the same reason, $U(na+mc)h\in\Cal Z$ for all $n,m\in\Bbb Z$.
Since $c/a\not\in\Bbb Q$, the subgroup $\{na+mc\mid n,m\in\Bbb Z\}$ is dense in $\Bbb R$.
Therefore, $U(t)h\in \Cal Z$ for all $t\in\Bbb R$.
It follows that $\Cal Z=\Cal H$.
\qed
\enddemo

\demo{Proof of Theorem 4.4}
It is enough to consider the Poissonian case  only
since the Gaussian case follows from it.

Suppose that we have a measure preserving  flow $T$ on a $\sigma$-finite infinite measure space $(X,\goth B,\mu)$ such that the following holds.
\roster
\item"(i)" $\exp(U_T(s))$ has a simple spectrum for each $s\in S$.
\item"(ii)"
For each finite sequence $s_1<s_2<\cdots<s_{k}$ of elements in  $S$
and each $1\le l_0\le k$,
there is a sequence of integers $t_j\to\infty$ such that
$$
\align
&U_T(t_js_l) \to \frac 1{2k}I\quad \text{if }1\le l\le k,\,l\ne l_0,\tag{4-3}\\
&U_T(t_{j}s_{l_0}) \to  \frac 1{2k}U_T(s_{l_0})\quad \text{if }j\to\infty\text{ and}\tag4-4\\
& U_T(t_jb)\to 0\quad\text{ for each positive }b\not\in S.\tag4-5
\endalign
$$
\endroster
We first show how to use this flow to prove the theorem and
 after that  we will explain how to construct such a flow.

Given a finite sequence of reals $0<z_1<\cdots< z_k$ and an integer vector $(n_1,\dots,n_k)\in\Bbb N^k$, we let
$O^{n_1,\dots,n_k}_{z_1,\dots,z_k}:=U_T(z_1)^{\odot n_1}\otimes\cdots \otimes U_T(z_k)^{\odot n_k}$.

{\bf Claim A.} For each  finite  sequence $s_1<\cdots< s_k$  of elements of $S$ and each $(n_1,\dots,n_k)\in\Bbb N^k$,
the operator
$O^{n_1,\dots,n_k}_{s_1,\dots,s_k}$
has a simple spectrum.

We verify this claim by induction in $k$.
If $k=1$ then the claim is true by (i).
Suppose it is true for some $k$.
Take a sequence $s_1<\cdots< s_{k+1}$ and $(n_1,\dots,n_{k+1})\in\Bbb N^{k+1}$.
By the inductive hypothesis, $O^{n_1,\dots,n_k}_{s_1,\dots,s_k}$
has a simple spectrum.
The operator $U_T(s_{k+1})^{\odot n_{k+1}}$ has a simple spectrum by (i).
Letting $l_0=k+1$ we deduce from \thetag{4-3} and \thetag{4-4} that
$$
\text{\rom{WCT}}(O^{n_1,\dots,n_k}_{s_1,\dots,s_k} \otimes U_T(s_{k+1})^{\odot n_{k+1}})\ni \frac 1{(2k+2)^{n_1+\cdots+n_{k+1}}}I\otimes U_T(s_{k+1})^{\odot n_{k+1}}.
$$
Now Lemma~4.6 yields that the operator $O^{n_1,\dots,n_{k+1}}_{s_1,\dots,s_{k+1}}=O^{n_1,\dots,n_k}_{s_1,\dots,s_k} \otimes U_T(s_{k+1})^{\odot n_{k+1}}$
has a simple spectrum.

{\bf Claim B.}
 Given two finite sequences $s_1<\cdots<s_k$ and $s_1'<\cdots< s_d'$ of elements from $S$ and two integer vectors
 $(n_1,\dots,n_k)\in\Bbb N^k$ and $(m_1,\dots,m_d)\in\Bbb N^d$,
if $\{s_1,\dots,s_k\}\ne\{s_1',\dots,s_d'\}$ then
$O^{n_1,\dots,n_k}_{s_1,\dots,s_k}$ is spectrally disjoint with
$O^{m_1,\dots,m_d}_{s_1',\dots,s_d'}$.

Without loss of generality we may assume that there is $1\le l_0\le k$ such that $s_{l_0}\not\in\{s'_1,\dots,s'_d\}$.
Then we deduce from \thetag{4-3} and \thetag{4-4} that there is a sequence of integers $t_j\to\infty$ such that
$$
\align
\big(O^{m_1,\dots,m_d}_{s_1',\dots,s_d'}
\big)^{t_j}   &\to \frac 1{(2r)^{m_1+\cdots+m_d}} I^{\otimes (m_1+\cdots+m_d)}\quad\text{and}\\
\big(O^{n_1,\dots,n_k}_{s_1,\dots,s_k}\big)^{t_j}
&\to \frac 1{(2r)^{n_1+\cdots+n_{k}}} I^{\otimes (n_1+\cdots+n_{l_0-1})}\otimes U_T(s_{l_0})^{\odot n_{l_0}}\otimes
I^{\otimes (n_{l_0+1}+\cdots+n_{k})}
\endalign
$$
as $j\to\infty$, where $r$ is the cardinality of the set $\{s_1,\dots,s_k,s_1',\dots,s_d'\}$.
Hence
the unitary operators
$O^{n_1,\dots,n_k}_{s_1,\dots,s_k}$ and  $O^{m_1,\dots,m_d}_{s_1',\dots,s_d'}$ are spectrally disjoint, as claimed.

{\bf Claim C.} Let $0<b\not\in S$.
Given two finite sequences $s_1<\cdots<s_k$ and $s_1'<\cdots< s_d'$ of elements from $S$ and two integer vectors
 $(n_1,\dots,n_k)\in\Bbb N^k$ and $(m_1,\dots,m_d)\in\Bbb N^d$,
the operators
$O^{n_1,\dots,n_k}_{s_1,\dots,s_k}$ and
$O^{m_1,\dots,m_d}_{bs_1',\dots,bs_d'}$
are spectrally disjoint.

Indeed, we deduce from \thetag{4-3} and \thetag{4-5} that there is a sequence of integers $t_j\to\infty$ such that
$$
\big(O^{n_1,\dots,n_k}_{s_1,\dots,s_k}\big)^{t_j}
\to \frac 1{(2r)^{n_1+\cdots+n_{k}}} I^{\otimes (n_1+\cdots+n_k)}\quad\text{and}\quad
\big(O^{m_1,\dots,m_d}_{ts_1',\dots,ts_d'}
\big)^{t_j}   \to 0
$$
as $j\to\infty$, where $r$ is the cardinality of the set $\{s_1,\dots,s_k,s_1',\dots,s_d'\}$.
The claim follows.

Now let $(Y,\nu)=(X,\mu)\times (S,\kappa)$, where $\kappa$ is the {\it counting} measure on $S$.
We define a flow $W=(W_t)_{t\in\Bbb R}$ on $(Y,\nu)$ by setting
$$
W_t(x,s):=(T_{st},s), \qquad x\in X, s\in S.
$$
Then $W$ preserves the $\sigma$-finite measure $\nu$.
This flow is not ergodic but every invariant subset is of either infinite or zero measure.
The Koopman representation of $\Bbb R$ associated with $W$ is $U_W=\bigoplus_{s\in S}U_T\circ c$.
Let $\widetilde W$ denote the Poisson suspension of $W$.
Since $U_{\widetilde W_1} =\exp(U_W(1))\ominus\Bbb C$ (see e.g., \cite{Ne}),
 we have
$$
\aligned
U_{\widetilde W_1}
&=\bigg(\bigotimes_{s\in S}\exp (U_{T}(s))\bigg)\ominus\Bbb C\\
&=\bigoplus_{k=1}^\infty\bigoplus_{s_1<\dots<s_k}
\bigoplus_{n_1=1}^\infty\cdots \bigoplus_{n_k=1}^\infty
O^{n_1,\dots,n_k}_{s_1,\dots,s_k},
\endaligned
\tag4-6
$$
where $s_1,\dots,s_k$ run over $S$.
It now follows from Claims A and B that the operator $U_{\widetilde W_1}$
has a simple spectrum.
Hence the flow $\widetilde W$ also has a simple spectrum.
It is obvious that $S\subset I(\widetilde W)$.

Now take a positive $b\not\in S$.
We are going to show that the Poisson flow $\widetilde W\circ b$ is spectrally disjoint with $\widetilde W$.
For that it is enough  to prove that the transformations $\widetilde W_b$
and $\widetilde W_1$ are spectrally disjoint.
We have
$$
\aligned
U_{\widetilde W_b} &
=\bigg(\bigotimes_{s\in S}\exp (U_{T}(bs))\bigg)\ominus\Bbb C\\
&=\bigoplus_{k=1}^\infty\bigoplus_{s_1<\dots<s_k}
\bigoplus_{n_1=1}^\infty\cdots \bigoplus_{n_k=1}^\infty
O^{n_1,\dots,n_k}_{bs_1,\dots,bs_k},
\endaligned
\tag4-7
$$
where $s_1,\dots,s_k$ run over $S$.
It remains to compare \thetag{4-6} and \thetag{4-7} and apply Claim~C.

To complete the proof of Theorem~4.4 we need  to construct the dynamical system $(X,\mu,T)$ satisfying (i) and (ii).
For that we use  the cutting-and-stacking inductive construction of rank-one flows.
The flow $T$ will be a rank-one flow associated with a sequence $(r_n,\sigma_n)_{n=1}^\infty$.
Thus our purpose is to define the sequence of cuts $r_n$ and spacer maps $\sigma_n:\{1,\dots,r_n\}\to\Bbb R_+$.
For that we  partition $\Bbb N$ into infinite subsets:
$$
\Bbb N=\bigg( \bigsqcup_{s\in S}\bigsqcup_{q\in\Bbb N}\bigsqcup_{i=1}^2\Cal L_{s,q}^i\bigg)\sqcup\bigg(
\bigsqcup_{k=1}^\infty
\bigsqcup_{s_1<\cdots<s_k}\bigsqcup_{l_0=1}^k
\Cal M_{s_1,\dots,s_k}^{l_0}\bigg),
$$
where $s_1,\dots,s_k$ run over $S$.

If for each $s\in S$ and $q\in\Bbb N$,
$$
\text{WCP}(U_T(s))\ni U_T(\sqrt 2 s)
\quad\text{and}\quad
\text{WCP}(U_T(s))\ni\frac 1q I+\frac{q-1}q U_T(s)
\tag4-8
$$
 then $U_T(s)$ has a simple spectrum by Lemma~4.7 and $\exp(U_T(s))$ has a simple spectrum by Lemma~4.5, i.e. (i) is satisfied.
To achieve this, we put
\roster
\item"---"  $r_n=n!$ and $\sigma_n(i)=\sqrt 2s$ for all $1\le i\le r_n$
for all $n\in\Cal L_{s,q}^1$,
\item"---"  $r_n=n!$ and $\sigma_n(i)=0$ if $1\le i<r_n/q$ and $\sigma_n(i)=s$ if $(1-q)r_n/q\le i\le r_n$ for all $n\in\Cal L_{s,q}^2$.
\endroster
A standard verification implies that
$U_T(-h_n)\to U_T(\sqrt 2 s)$ if $\Cal L_{s,q}^1\ni n\to\infty$ and
$U_T(-h_n)\to \frac 1q I+\frac{q-1}q U_T(s)$ if $\Cal L_{s,q}^2\ni n\to\infty$, where $h_n$ as usual denotes the height of the $n$-th tower.
We note that though $T$ has not yet been defined entirely, these limits are well defined because  they do  not depend on the choice of $r_n,\sigma_n$
when $n\not\in\bigsqcup_{s\in S}\bigsqcup_{q\in \Bbb N}(\Cal L_{s,q}^1\sqcup\Cal L_{s,q}^2)$.
Thus, \thetag{4-8}, and hence (i), is satisfied.

To realize (ii) we fix a
finite sequence $s_1<s_2<\cdots<s_{k}$ of elements in  $S$
and $1\le l_0\le k$.
Enumerate the elements of $\Cal M_{s_1,\dots,s_k}^{l_0}$ in ascending  order: $n_1<n_2<\cdots$.
We now let $r_{n_j}:=2k$ for all $j$.
Instead of writing  precise formulas for the spacer maps $\sigma_{n_j}:\{1,\dots,2k\}\to\Bbb R_+$ we illustrate
the idea of the construction with the following picture of the $(n_j+1)$-th tower in this subsequence (see Figure 4.1).
To be specific, we choose $k=3$ and $l_0=2$.
Since the tower is very ``high'', we place it horizontally.

\midinsert
\beginpicture
\put {} at 0 0
\put {} at 0 20
\plot 1 0   360 0  360 20 0 20  /
\setplotsymbol({\llarge .})
\plot 0 1    0 19 /
\plot 8 1    8 19 /
\plot 40 1   40 19 /
\plot 55 1   55 19 /
\plot 125 1   125 19 /
\plot 200 1   200 19 /
\setdots <1pt>
\setplotsymbol({\tiny .})
\plot 0 0 0 -8 /
\plot 8 0 8 -8 /
\plot 40 0 40 -8 /
\plot 55 0 55 -8 /
\plot 125 0 125 -8 /
\plot 200 0 200 -8 /
\setsolid
\arrow <3pt> [0.5, 1] from 0 -8 to 8 -8
\arrow <3pt> [0.5, 1] from 8 -8 to 0 -8
\arrow <3pt> [0.5, 1] from 40 -8 to 55 -8
\arrow <3pt> [0.5, 1] from 55 -8 to 40 -8
\arrow <3pt> [0.5, 1] from 125 -8 to 200 -8
\arrow <3pt> [0.5, 1] from 200 -8 to 125 -8
\put {$_{t_js_1}$} at 5 -14
\put {$_{t_js_2-s_{l_0}}$} at 47 -14
\put {$_{t_js_3}$} at 163 -14
\endpicture
\botcaption{Figure 4.1}
$(n_j+1)$-th tower.
\endcaption
\endinsert
\noindent
The black stripes here are the copy of the $n_j$-th tower.
They are very ``thin'' because we choose the parameter $t_j\in\Bbb N$ very large.
It is easy to see that \thetag{4-3} and \thetag{4-4} hold.
Denote by $a_{j,i}$ the distances between  the $2i$-th and $(2i+1)$-th copies of the $n_j$-th tower in the $(n_j+1)$-th tower, $j=1,2$.
Let $a_{j,3}$ be the distance between the $6$-th copy and the top of the
$(n_j+1)$-th tower.
We arrange the spacers in the $(n_j+1)$-th tower in such a way
that $t_js_1 \ll a_{j,1}\ll a_{j,2}\ll a_{j,3}$, where the sign ``$\ll$'' means {\it grows much faster} as $j\to \infty$.
Then \thetag{4-5} follows.
\qed
\enddemo

Using Theorem~4.4, we can  sharpen Lemma~3.1.

\proclaim{Corollary 4.8}
Let $\Cal P_C^0$ denote the set of continuous, fully supported measures
$\sigma$ on $\Bbb R$ such that $(\sigma^p)_t\perp\sigma^q$ for all $t>0$,
$t\ne 1$, and $p,q\in\Bbb N$, where the upper indices  $p,q$ denote the convolution powers.
Then $\Cal P_C^0$ is a dense $G_\delta$ in $\Cal P$.
\endproclaim

\demo{Proof}
We first recall a well known fact that the fully supported continuous measures on $\Bbb R$ form a $G_\delta$ subset in $\Cal P$ (see, e.g. \cite{Na}).
Since the maps $\Cal P\ni\sigma\mapsto\sigma^p\in\Cal P$ are continuous for all $p\in\Bbb N$,
we can argue as in the proof of lemma~3.1 to show that $\Cal P_C^0$ is a $G_\delta$.
It follows from the proof of Theorem~4.4 that a  measure  of maximal spectral type of $U_W$ belongs to $\Cal P_C^0$.
It remains to note that the equivalence class of every fully supported non-atomic measure is dense in $\Cal P$.
\qed
\enddemo

\subhead III. Asymmetries in flows
 \endsubhead
If $F$ is a Gaussian flow (or a single transformation)  then $F$ is isomorphic to its inverse and hence $-1\in I(F)$.
In contrast to this, there exist {\it asymmetric} Poisson suspensions, i.e. Poisson suspension flows $\widetilde T$ such that $-1\not \in I(\widetilde T)$.
To construct such a $\widetilde T$, we use the ideas from \cite{Ry5} and \cite{Ro}.

\proclaim{Proposition 4.9}
Let $T$ be an  infinite measure preserving rank-one flow associated with
a sequence $(r_n,s_n)_{n=1}^\infty$ of cutting-and-stacking parameters.
If for some sequence $l_i\to\infty$,
$$
\aligned
&r_{l_i}=5\text{ and } s_{l_i}(0)=0,\ s_{l_i}(1)=s_{l_i}(2)=1, \ s_{l_i}(3)=s_{l_i}(4)=2,\\
& r_{l_i+1}\to\infty \text{ and }s_{l_i+1}\equiv 0,
\endaligned
\tag4-9
$$
then the transformation $T_1$ is not isomorphic to $T_{-1}$.
If, in addition, $T$ has a simple spectrum then the Poisson suspension flow $\widetilde T$ is asymmetric.
\endproclaim
\demo{Proof}
Let $n_i:=h_{l_i}+1$.
We recall that $h_{l_i}$ is the height of the $l_i$-th tower in the inductive construction of $T$.
Repeating the argument from \cite{Ry5} almost verbally,
we obtain that
$$
\align
&\lim_{i\to\infty}\mu(A\cap T_{n(i)}A\cap T_{3n(i)}A)\ge 0.2\mu(A)\quad \text{for all subsets $A$, $\mu(A)<\infty$, and }\\
&\lim_{i\to\infty}\mu(A'\cap T_{-n(i)}A'\cap T_{-3n(i)}A')=0\quad \text{for some subset $A'$, $0<\mu(A')<\infty$}.
\endalign
$$
This implies that $-1\not\in E(T)$, as desired.

 If  $T$ has a simple spectrum and $\widetilde T_1$ is isomorphic to $\widetilde T_{-1}$
then it follows from \cite{Ro, Proposition~5.2} that   $T_1$ is isomorphic to $T_{-1}$, a contradiction.
Thus $\widetilde T$ is asymmetric.
\qed
\enddemo

 We can now refine the first claim of Theorem~4.4 in the following way.

\proclaim{Theorem 4.10}
 Let $S$ be a countable subgroup of $\Bbb R^*$.
There is a weakly mixing Poisson suspension  flow $\widetilde W$ with a simple spectrum such that $I(\widetilde W)=S$ and $\sigma_{\widetilde W}\perp(\sigma_{\widetilde W})_t$ for each positive $t\not\in S$.
\endproclaim

\demo{Proof}
 We consider separately two cases.

{\bf Case 1.} Let $S\subset \Bbb R^*_+$.
Then we construct $T$ as in the proof of Theorem~4.4 but add an extra condition \thetag{4-9} on the sequence of cutting-and-stacking parameters.
Define the flow $W$ by $T$ and $S$ in the same way as in the proof of Theorem~4.4.
Then the only what we need to show is that $-1\not\in I(\widetilde W)$.
Indeed, if $-1\in I(\widetilde W)$ then by \cite{Ro} there $RW_1R^{-1}=W_{-1}$ for a measure preserving transformation $R$.
We note that the transformation $W_1$ is not ergodic.
The maximal spectral types of all its ergodic components are pairwise  orthogonal.
It follows that $R$ preserves every ergodic component of $W_1$.
Hence every ergodic component of $W_1$ is conjugate to its inverse.
In particular, $T_1$ is conjugate to $T_{-1}$.
However this contradicts to Proposition~4.9.

{\bf Case 2.} Let $S=S'\sqcup(-S')$, where $S'$ is a countable subgroup of $\Bbb R^*_+$.
Then we use a  {\it symmetrization trick}.
First, we construct an infinite measure preserving flow $T$  exactly as in the proof of Theorem~4.4 but with $S'$ instead of $S$.
Since $T$ is a rank-one flow, it is associated with a sequence of
 cutting-and-stacking parameters $(r_n,s_n)_{n=1}^\infty$.
We now consider a new sequence $(r_n',s_n')_{n=1}^\infty$,
where $r_n':=2r_n-1$ and $s_n'$ is a map from the set $\{-1,0,\dots, r_n'-1\}$ to $\Bbb R_+$ given by
$$
s_n'(i)=\cases
s_n(r_n-2-i)&\text{if }-1\le i\le r_n-2,\\
s_n(i-r_n+1)&\text{if }r_n-1\le i\le r_n'-1.
\endcases
$$
Let $T'$ denote the rank-one flow associated with $(r_n',s_n')_{n=1}^\infty$.
It is defined in a usual way with one exception.
For each $n$, we  have an {\it extra} value of $s_n'$ at the point $-1$.
This means that  we first construct a preliminary $(n+1)$-th tower using the $n$-th tower, $r_n'$ and $s_n'(i)$, $0\le i<r_n'$.
 Then  we enlarge it by adding an additional {\it spacer} rectangular of height $s_n'(-1)$ underneath of it.
This new tower is the $(n+1)$-tower of the inductive construction for $T'$ (see Figure 4.2, where $r_n=2$).

\midinsert
\beginpicture
\put {} at 0 0
\put {} at 0 25
\setplotsymbol({\tiny .})

\plot 0 0   260 0  260 25 0 25  /
\plot 0 0    0 25 /
\plot 10 0    10 25 /
%\setplotsymbol({\llarge $\cdot$})
\plot 80 0   80 25 /
\plot 95 0   95 25 /
\plot 165 0   165 25 /
\plot 180 0    180 25 /
\plot 250 0    250 25 /
\plot 260 0    260 25 /
\setdots <1pt>
\plot 0 0 0 -8 /
\plot 10 0 10 -8 /
\plot 95 0 95 -8 /
\plot 80 0 80 -8 /
\plot 130 -5 130 30 /
\plot 250 0 250 -8 /
\plot 260 0 260 -8 /
\plot 165 0 165 -8 /
\plot 180 0 180 -8 /

\setsolid
\arrow <3pt> [0.5, 1] from 0 -8 to 10 -8
\arrow <3pt> [0.5, 1] from 10 -8 to 0 -8
\arrow <3pt> [0.5, 1] from 80 -8 to 95 -8
\arrow <3pt> [0.5, 1] from 95 -8 to 80 -8
\arrow <3pt> [0.5, 1] from 165 -8 to 180 -8
\arrow <3pt> [0.5, 1] from 180 -8 to 165 -8
\arrow <3pt> [0.5, 1] from 260 -8 to 250 -8
\arrow <3pt> [0.5, 1] from 250 -8 to 260 -8
\put {$_{s_n(1)}$} at 5 -14
\put {$_{s_n(0)}$} at 88 -14
\put {$_{s_n(0)}$} at 173 -14
\put {$_{s_n(1)}$} at 256 -14
\endpicture
\botcaption{Figure 4.2}
$(n+1)$-th symmetrized tower.
\endcaption
\endinsert
\noindent
It is straightforward to verify that \thetag{4-8}, \thetag{4-5} and analogues of \thetag{4-3} and \thetag{4-4} with coefficients $\frac 2{4k+1}$ instead of $\frac 1{2k}$ hold for $T'$ (with $S'$ instead of $S$).
Therefore, if we define a flow $W$ in the same way as in the proof of Theorem~4.4 but with $S'$ instead of $S$ then the Poisson suspension $\widetilde W$ of $W$  has a simple spectrum, $I(\widetilde W)\cap\Bbb R^*_+=S'$ and $\sigma_{\widetilde W}\perp(\sigma_{\widetilde W})_t$ for all positive $t\not\in S'$.
Next, we note that $T$ is conjugate to its inverse.
The corresponding conjugation $R$ can be defined inductively in the following way.
Let $X_n$ be an $n$-th tower and let $A$ be a  subset in $X_n$ such that if $x,y\in A$ are on the same height in $X_n$ and $x\in A$ then $y\in A$.
We denote by $A^*$ the subset of $X_n$ which is symmetric to $A$ with respect to the horizontal line passing through   the middle of the tower.
Then we set $RA:=A^*$.
Passing to the limit when $n\to\infty$ we obtain
a well-defined invertible measure preserving transformation $R$ of $X$ and $RT_tR^{-1}=T_{-t}$ for all $t\in\Bbb R$.
We recall that the flow $W$ is defined on the space $(X\times S',\mu\times\kappa)$.
We now define a transformation $Q$ of this space by setting
$Q(x,s)=(Rx,s)$.
Then $QW_tQ^{-1}=W_{-t}$ for all $t$.
This implies that $-1\in I(\widetilde W)$.
Hence $I(\widetilde W)=S$, as desired. \qed
\enddemo

\remark{Remark \rom{4.11}} If we drop the condition in Theorem~4.10 that $\widetilde W$ has a simple spectrum then the proof of Case 2 simplifies.
Indeed, we do not need the symmetrization trick.
It is enough to argue as in the proof of Theorem~4.4.
The flow $W$ will act on the space $X\times S$.
We then define the transformation $Q$ conjugating $W$ with $W\circ(-1)$ by $Q(x,s)=(x,-s)$.
\endremark

\remark{Remark \rom{4.12}}
It is possible to  strengthen Theorems~4.4 and 4.10 by constructing a {\it mixing}
Poisson suspension (and a {\it mixing} Gaussian) flow
$\widetilde W$ satisfying the conditions of those theorems.
For that one should apply  the technique of {\it forcing  of mixing} developed in  our previous paper \cite{DaR} (see also \cite{Ry6}).
\endremark

\enddocument